\documentclass[12pt]{amsart}
\usepackage{amsfonts, amssymb}
\usepackage{amsmath}
\usepackage{amsthm}
\usepackage{cite}
\numberwithin{equation}{section}

\newtheorem{theorem}{Theorem}[section]
\newtheorem{lemma}[theorem]{Lemma}

\newtheorem{corollary}[theorem]{Corollary}

\begin{document}

\title[Oscillation of differential equations]{On the oscillation of certain second-order linear differential equations}

\author{Yueyang Zhang}
\address{School of Mathematics and Physics, University of Science and Technology Beijing, No.~30 Xueyuan Road, Haidian, Beijing, 100083, P.R. China}
\email{zyynszbd@163.com}





\date{\today}

\commby{}

\begin{abstract}
This paper consists of three parts: First, letting $b_1(z)$, $b_2(z)$, $p_1(z)$ and $p_2(z)$ be nonzero polynomials such that $p_1(z)$ and $p_2(z)$ have the same degree $k\geq 1$ and distinct leading coefficients $1$ and $\alpha$, respectively, we solve entire solutions of the Tumura--Clunie type differential equation $f^{n}+P(z,f)=b_1(z)e^{p_1(z)}+b_2(z)e^{p_2(z)}$, where $n\geq 2$ is an integer, $P(z,f)$ is a differential polynomial in $f$ of degree $\leq n-1$ with coefficients having polynomial growth. Second, we study the oscillation of the second-order differential equation $f''-[b_1(z)e^{p_1(z)}+b_2(z)e^{p_2(z)}]f=0$ and prove that $\alpha=[2(m+1)-1]/[2(m+1)]$ for some integer $m\geq 0$ if this equation admits a nontrivial solution such that $\lambda(f)<\infty$. This partially answers a question of Ishizaki. Finally, letting $b_2\not=0$ and $b_3$ be constants and $l$ and $s$ be relatively prime integers such that $l> s\geq 1$, we prove that $l=2$ if the equation $f''-(e^{lz}+b_2e^{sz}+b_3)f=0$ admits two linearly independent solutions $f_1$ and $f_2$ such that $\max\{\lambda(f_1),\lambda(f_2)\}<\infty$. In particular, we precisely characterize all solutions such that $\lambda(f)<\infty$ when $l=2$ and $l=4$.
\end{abstract}

\maketitle


\section{Introduction}\label{se1: introduction}

In the last several decades, the growth and value distribution of meromorphic solutions of complex differential equations have attracted much interest; see \cite{Laine1993} and references therein. One of the main tools in this subject is Nevanlinna theory; see, e.g., \cite{Hayman1964Meromorphic,Laine1993} for the standard notation and basic results of Nevanlinna theory. Bank and Laine~\cite{Banklaine1982,Banklaine1982-2} initiated the study on the oscillation of the second-order linear differential equation
\begin{equation}\label{bank-laine0}
f''+A(z)f=0,
\end{equation}
where $A(z)$ is an entire function. It is well-known that all solutions of equation \eqref{bank-laine0} are entire. For an entire function $f$, denote by $\sigma(f)$ the \emph{order} of $f$ which is defined as
\begin{equation*}
\sigma(f)=\limsup_{r\to\infty}\frac{\log T(r,f)}{\log r}=\limsup_{r\to\infty}\frac{\log\log M(r,f)}{\log r},
\end{equation*}
where $M(r,f)$ is the maximum modulus of $f$ on the circle $|z|=r$. When $A$ is transcendental, an application of the lemma on the logarithmic derivative easily yields that all nontrivial solutions of \eqref{bank-laine0} satisfy $\sigma(f)=\infty$. Denote by $\lambda(f)$ the \emph{exponent of convergence of zeros} of $f$ which is defined as
\begin{equation*}
\lambda(f)=\limsup_{r\to\infty}\frac{\log n(r,f)}{\log r},
\end{equation*}
where $n(r,f)$ denotes the number of zeros of $f$ in the disc $\{z: |z|<r\}$. Concerning the zero distribution of solutions of equation \eqref{bank-laine0}, Bank and Laine \cite{Banklaine1982,Banklaine1982-2} proved: Let $f_1$ and $f_2$ be two linearly independent solutions of \eqref{bank-laine0}. If $\sigma(A)$ is not an integer, then $\max\{\lambda(f_1),\lambda(f_2)\}\geq \sigma(A)$; if $\sigma(A)<1/2$, then $\max\{\lambda(f_1),\lambda(f_2)\}=\infty$. Later, Shen~\cite{Shen1985} and Rossi~\cite{Rossi1986} relaxed the condition $\sigma(A)<1/2$ to the case $\sigma(A)=1/2$. Based on these results, Bank and Laine conjectured that $\max\{\lambda(f_1),\lambda(f_2)\}=\infty$ whenever $\sigma(A)$ is not an integer. This conjecture is known as \emph{the Bank--Laine conjecture} and has attracted much interest; see the surveys \cite{Gundersen2014,lainetohge2008} and references therein. Recently, this conjecture was disproved by Bergweiler and Eremenko \cite{Bergweilereremenko2017,Bergweilereremenko2019}. They constructed counterexamples for the coefficient $A$ such that $\sigma(A)$ is not an integer and equation \eqref{bank-laine0} admits two linearly independent solutions such that $\max\{\lambda(f_1),\lambda(f_2)\}<\infty$. In particular, one of the solutions is free of zeros. In their constructions, they used the solutions of \eqref{bank-laine0} with $A$ being a polynomial of $e^z$ of degree~2, namely $A(z)=a_1e^{2z}+a_2e^z+a_3$ with certain coefficients $a_1$, $a_2$ and $a_3$.

On the other hand, it is natural to give explicit solutions of \eqref{bank-laine0} such that $\lambda(f)<\infty$ when $A$ is a periodic entire function of the form
\begin{equation}\label{bank-laine01}
A(z)=B(e^z), \quad B(\zeta)=b_{-k}\zeta^{-k}+\cdots+b_0+\cdots+b_l\zeta^{l}, \quad b_{-k}b_l\not=0.
\end{equation}
For such solutions, a remarkable result in \cite{Banklaine1983-1,Chiang2000} states that there exist complex constants $c$, $c_j$ and a polynomial $P(z)$ with simple roots only such that if $l$ is an odd positive integer, then
\begin{equation}\label{bank-laine01hf1}
f=P(e^{z/2})\exp\left(\sum_{j=0}^{l}c_je^{(l-j)z/2}+cz\right),
\end{equation}
where $c_j=0$ whenever $j$ is even; while if $l$ is an even positive integer, then
\begin{equation}\label{bank-laine01hf2}
f=P(e^{z})\exp\left(\sum_{j=0}^{l/2}c_je^{(l/2-j)z}+cz\right).
\end{equation}
However, it seems difficult to determine explicitly $c_j$ and also the polynomial $P(z)$ in the above two expressions and, until now, they are only known in some special cases. For example, Bank and Laine~\cite{Banklaine1983-1} gave a precise characterization of all nontrivial solutions such that $\lambda(f)<\infty$ of \eqref{bank-laine0} when $A(z)=e^{z}-b$ for some constant $b$; see also \cite[Theorem~5.22]{Laine1993}. Bank and Laine~\cite{Banklaine1983-1} also characterized entire solutions such that $\lambda(f)<\infty$ of equation \eqref{bank-laine0} when $A(z)=-(1/4)e^{-2z}+(1/2)e^{-z}+b$ for some constant $b$. For these two coefficients, Chiang and Ismail~\cite{Chiang2006} expressed all solutions of \eqref{bank-laine0} in terms of some special functions and give a complete characterization of the zero distribution of these solutions.

In \cite{Bank1993-1}, Bank developed a method to find entire solutions such that $\lambda(f)<\infty$ of equation \eqref{bank-laine0}, but the manipulation of this method seems complicated. One of the main purposes of this paper is to give a more precise description of the oscillation of equation \eqref{bank-laine0} when $A(z)$ contains two exponential terms, i.e.,
\begin{equation}\label{bank-laine011}
A(z)=B(e^z), \quad B(\zeta)=b_{-k}\zeta^{-k}+b_0+b_l\zeta^{l}, \quad b_{-k}b_l\not=0,
\end{equation}
or
\begin{equation}\label{bank-laine012}
A(z)=B(e^z), \quad B(\zeta)=b_0+b_{s}\zeta^{s}+b_l\zeta^{l}, \quad b_sb_l\not=0.
\end{equation}
In particular, this provides a different approach from that in \cite{Chiang2006} and also leads to a complete characterization of all solutions such that $\lambda(f)<\infty$ of \eqref{bank-laine0} when $A(z)$ is an arbitrary polynomial in $e^z$ of degree~$2$; see Theorem~\ref{maintheorem6} in section~\ref{se4: equation with periodic coefficients}. This work is a continuation of~\cite{zhang2021-1}, where the present author found all nontrivial solutions such that $\lambda(f)<k$ of the differential equation
\begin{equation}\label{bank-laine}
f''-\left[b_1(z)e^{p_1(z)}+b_2(z)e^{p_2(z)}+b_3(z)\right]f=0,
\end{equation}
where $b_1(z)$, $b_2(z)$ and $b_3(z)$ are three polynomials such that $b_1(z)b_2(z)\not\equiv0$ and $p_1(z)$ and $p_2(z)$ are two polynomials of the same degree~$k\geq 1$ with distinct leading coefficients $1$ and $\alpha$, respectively.

\begin{theorem}[see \cite{zhang2021-1}]\label{maintheorem1}
Let $b_1$, $b_2$ and $b_3$ be polynomials such that $b_1b_2\not\equiv0$ and $p_1$, $p_2$ be two polynomials of degree~$k\geq 1$ with distinct leading coefficients $1$ and $\alpha$, respectively, and~$p_1(0)=p_2(0)=0$. Suppose that \eqref{bank-laine} admits a nontrivial solution such that $\lambda(f)<k$. Then $\alpha=1/2$ or $\alpha=3/4$. Moreover,
\begin{itemize}
\item [(1)]
if $\alpha=1/2$, then $p_2=p_1/2$, $f=\kappa e^{h}$, where $\kappa$ is a polynomial with simple roots only and $h$ satisfies $h'=\gamma_1e^{p_1/2}+\gamma$ with $\gamma_1$ and $\gamma$ being two polynomials such that $\gamma_1^2=b_1$, $2\gamma_1\gamma+\gamma_1'+\gamma_1p_1'/2+2\kappa'/\kappa\gamma_1=b_2$ and $\gamma^2+\gamma'+2\gamma\kappa'/\kappa+\kappa''/\kappa=b_3$;

\item [(2)]
if $\alpha=3/4$, then $p_1=z$, $p_2=3z/4$ and $f=e^{h}$, where $h$ satisfies $h'=-4c^2e^{z/2}+ce^{z/4}-1/8$ and $A=-(16c^2e^{z}-8c^3e^{3z/4}+1/64)$, where $c$ is a nonzero constant.
\end{itemize}
\end{theorem}

The proof of Theorem~\ref{maintheorem1} is based on a development of the Tumura--Clunie method; see \cite[Chapter~4]{Hayman1964Meromorphic}. Define a \emph{differential polynomial} $P(z,g)$ in $g$ to be a finite sum of monomials in $g$ and its derivatives of the form $P(z,g)=\sum_{l=1}^{m}a_{l}g^{n_{l0}}(g')^{n_{l1}}\cdots(g^{(s)})^{n_{ls}}$,
where $n_{l0},\cdots,n_{ls}\in \mathbb{N}$ and the coefficients $a_l$ are meromorphic functions of order less than~$\sigma(g)$. Define the \emph{degree} of $P(z,g)$ to be the greatest integer of $d_l:=\sum_{t=0}^sn_{lt}$, $l=1,\cdots,m$, and denote it by $\deg_g(P(z,g))$. Consider the equation
\begin{equation}\label{EQ1 intr}
g^n+P(z,g)=b_1e^{p_1}+b_2e^{p_2},
\end{equation}
where $n\geq 2$ and $P(z,g)$ is a differential polynomial in $g$ of degree $\leq n-1$ with meromorphic functions of order less than~$k$ as coefficients. If equation~\eqref{bank-laine} admits an entire solution such that $\lambda(f)<k$, then equation \eqref{bank-laine} reduces to an equation of the form in \eqref{EQ1 intr} with $n=2$. It is shown in~\cite[Theorem~2.1]{zhang2021-1} that if equation \eqref{EQ1 intr} admits an entire solution, then either $\alpha=-1$ or $\alpha$ is positive rational number and in either case $g$ is a linear combination of certain exponential functions plus some function of order less than~$k$. However, to solve entire solutions of \eqref{bank-laine} such that $\lambda(f)<\infty$, \cite[Theorem~2.1]{zhang2021-1} fails to work since in this case the coefficients of $P(z,g)$ shall contain some logarithmic derivatives which have order no less than~$k$.

The remainder of this paper is organized in the following way. Denote by $\mathcal{R}$ the set of rational functions and by $\mathcal{L}$ the set of functions $a(z)$ such that $a(z)=h^{(l)}(z)/h(z)$, $l\geq 1$, for some meromorphic function $h(z)$ of finite order, respectively. In section~\ref{se2: Tumura--Clunie differential equations}, we further develop the Tumura--Clunie method by solving entire solutions of equation \eqref{EQ1 intr}, where $P(z,g)$ is now a differential polynomial in $g$ with coefficients that are combinations of functions in the set $\mathcal{S}=\mathcal{R}\cup\mathcal{L}$. For equation \eqref{EQ1 intr} with such coefficients, we can also write the entire solution as a linear combination of exponential functions with certain constant coefficients, but unlike in \cite[Theorem~2.1]{zhang2021-1}, it is impossible to determine whether $\alpha$ is a rational number; see Theorem~\ref{maintheorem3}. In section~\ref{se3: an oscillation question of Ishizaki}, we apply our results on equation \eqref{EQ1 intr} to study the oscillation of equation \eqref{bank-laine} and prove that $\alpha=[2(m+1)-1]/[2(m+1)$] for some integer $m\geq 0$ provided that equation \eqref{bank-laine} with $b_3\equiv0$ admits a nontrivial solution such that $\lambda(f)<\infty$; see Theorem~\ref{maintheorem4}. This gives a partial answer to a question of Ishizaki~\cite{Ishizaki19970}. In section~\ref{se4: equation with periodic coefficients}, we consider the equation $f''-(b_1e^{lz}+b_2e^{sz}+b_3)f=0$, where $l,s$ are relatively prime integers such that $l>s\geq 1$ and $b_i$ are constants such that $b_1b_2\not=0$. We prove that $l=2$ if this equation admits two linearly independent solutions $f_1$ and $f_2$ such that $\max\{\lambda(f_1),\lambda(f_2)\}<\infty$. In particular, when $l=2$ or $l=4$, we determine the polynomial $P(z)$ and the coefficients $c_j$ and $c$ in \eqref{bank-laine01hf2} precisely. Finally, in section~\ref{se5: concluding remarks}, we give some remarks on our results.

\section{Tumura--Clunie differential equations}\label{se2: Tumura--Clunie differential equations}

Let $b_1(z)$ and $b_2(z)$ be two nonzero polynomials and $p_1(z)$ and $p_2(z)$ be two polynomials of the same degree~$k\geq 1$ with distinct leading coefficients $1$ and $\alpha$, respectively, and $p_1(0)=p_2(0)=0$. Without loss of generality, we may suppose that $0<|\alpha|\leq 1$. In this section, we solve entire solutions of the differential equation
\begin{equation}\label{EQ1}
f^n+P(z,f)=b_1e^{p_1}+b_2e^{p_2},
\end{equation}
where $n\geq 2$ and $P(z,f)$ is a differential polynomial in $f$ of degree $\leq n-1$ with coefficients being combinations of functions in $\mathcal{S}$. In the following, a differential polynomial in $f$ will always have coefficients which are combinations of functions in $\mathcal{S}$ and thus we will omit mentioning this from now on.

To state our results, we first set up some notation: Let $p(z)$ be a polynomial of degree $k\geq 1$. We write $p(z)=(a+ib)z^{k}+q(z)$, where $a,b$ are real and $a+ib\not=0$ and $q(z)$ is a polynomial of degree at most $k-1$. Denote
\begin{equation}\label{expon1}
\delta(p,\theta)=a\cos k\theta-b\sin k\theta, \quad \theta\in[0,2\pi).
\end{equation}
Then on the ray $z=re^{i\theta}$, $r\geq 0$, from \cite{Banklangley1987} (or \cite[Lemma~5.14]{Laine1993}) we know that:
\begin{itemize}
  \item [1.] if $\delta(p,\theta)>0$, then there exists an $r_0=r_0(\theta)$ such that $\log |e^{p(z)}|$ is increasing on $[r_0,\infty)$ and $|e^{p(z)}|\geq e^{\delta(p,\theta)r^{n}/2}$ there;
  \item [2.] if $\delta(p,\theta)<0$, then there exists an $r_0=r_0(\theta)$ such that $\log |e^{p(z)}|$ is decreasing on $[r_0,\infty)$ and $|e^{p(z)}|\leq e^{\delta(p,\theta)r^{n}/2}$ there.
\end{itemize}
Let $\theta_1, \theta_2, \cdots, \theta_{2k}\in[0,2\pi)$ be such that $\delta(p,\theta_j)=0$, $j=1,2,\cdots,2k$. We may suppose that $\theta_1< \pi$ and $\theta_j=\theta_1+(j-1)\pi/k$. Denoting $\theta_{2k+1}=\theta_1+2\pi$, then $\theta_1$, $\theta_2$, $\cdots$, $\theta_{2k}$ divides the complex plane $\mathbb{C}$ into $2k$ sectors $S_j$, namely
\begin{equation}\label{expon10}
S_j=\left\{re^{i\theta}:\, 0\leq r<\infty, \quad \theta_j< \theta<\theta_{j+1}\right\}, \quad j=1,2,\cdots,2k.
\end{equation}
Throughout this paper, we let $\epsilon>0$ be an arbitrary constant. We also denote
\begin{equation}\label{expon10 fujia}
S_{j,\epsilon}=\left\{re^{i\theta}:\, 0\leq r<\infty, \quad \theta_j+\epsilon<\theta<\theta_{j+1}-\epsilon\right\}, \quad j=1,2,\cdots,2k.
\end{equation}
Denote by $\overline{S}_j$ and $\overline{S}_{j,\epsilon}$ the closure of $S_j$ and $S_{j,\epsilon}$, respectively. For $p_1$ in \eqref{EQ1}, we choose $\theta_1=-\pi/(2k)$ and thus $\delta(p_1,\theta)>0$ in the sectors $S_j$ when $j$ is odd, and $\delta(p_1,\theta)<0$ in the sectors $S_j$ when $j$ is even. Denote by $J_1$ and $J_2$ the subsets of odd and even integers in the set $J=\{1,2\cdots,2k\}$, respectively, i.e., $J_1=\{1,3,\cdots,2k-1\}$ and $J_2=\{2,4,\cdots,2k\}$. We prove the following

\begin{theorem}\label{maintheorem3}
Let $n\geq2$ be an integer and $P(z,f)$ be a differential polynomial in $f$ of degree $\leq{n-1}$. Suppose that \eqref{EQ1} admits an entire solution~$f$. Then $\alpha$ is real. Moreover,
\begin{itemize}
  \item [(1)]
if $-1\leq \alpha<0$, then $f=\gamma_1e^{p_1/n}+\gamma_2e^{p_2/n}+\eta$, where $\gamma_1$, $\gamma_2$ are two polynomials such that $\gamma_1^n=b_1$, $\gamma_2^n=b_2$ and $\eta$ is an entire function such that $\eta=(\mu_{1,j}-1)\gamma_1e^{p_1/n}+(\mu_{2,j}-1)\gamma_2e^{p_2/n}+\eta_j$, where $\mu_{1,j}$ and $\mu_{2,j}$ are the $n$-th roots of~$1$ such that $\mu_{1,j}=1$ when $j\in\{1\}\cup J_2$ and $\mu_{2,j}=1$ when $j\in\{2\}\cup J_1$, and there is an integer $N$ such that $|\eta_j|=O(r^N)$ uniformly in $\overline{S}_{j,\epsilon}$; in particular, when $k=1$, $\eta$ is a polynomial;

  \item [(2)]

if $0<\alpha<1$, letting $m$ be the smallest integer such that $\alpha\leq[(m+1)n-1]/[(m+1)n]$, then $f=\gamma_1\sum_{j=0}^mc_j(b_2/b_1)^je^{[jn(\alpha-1)+1]p_1/n}+\eta$, where $\gamma_1$ is a polynomial such that $\gamma_1^n=b_1$ and $c_0$, $\cdots$, $c_m$ are constants such that $c_0^{n}=1$ when $m=0$, and $c_0^{n}=nc_0^{n-1}c_1=1$ when $m=1$, and $c_0^{n}=nc_0^{n-1}c_1=1$ and $\sum_{\substack{j_0+\cdots+j_m=n,\\j_1+\cdots+mj_m=k_0}}\frac{n!}{j_0!j_1!\cdots j_m!}c_0^{j_0}c_1^{j_1}\cdots c_m^{j_m}=0$, $k_0=2,\cdots,m$, when $m\geq 2$, and $\eta$ is a meromorphic function with at most finitely many poles such that $\eta=\gamma_1\sum_{l=0}^m(\mu_j-1)c_{j}(b_2/b_1)^je^{[jn(\alpha-1)+1]p_1/n}+\eta_{j}$, where $\mu_j$ are the $n$-th roots of~$1$ such that $\mu_j=1$ when $j\in\{1\}\cup J_2$, and there is an integer $N$ such that $|\eta_j|=O(r^N)$ uniformly in $\overline{S}_{j,\epsilon}$; moreover, we have $p_2=\alpha p_1$ when $m\geq 1$; in particular, when $k=1$, $\eta$ is a rational function.
\end{itemize}
\end{theorem}

In Theorem~\ref{maintheorem3}, if all coefficients of the monomials in $P(z,f)$ of degree $n-1$ are rational functions, then we may use the method in the proof of \cite[Theorem~2.1]{zhang2021-1} to show that $\eta$ is a polynomial or a rational function. We also remark that, by using the method in the proof of Theorem~\ref{maintheorem3} for the case $-1\leq\alpha<0$ together with the method in \cite{zhang2021}, we may extend \cite[Theorem~2.1]{zhang2021-1} to the case $P(z,f)$ is a delay--differential polynomial in $f$ with meromorphic functions of order less than~$k$ as coefficients; see \cite{zhang2021} for the definition of a  delay--differential polynomial.

As in the proof of~Theorems~\cite[Theorem~1.1]{zhang2021} and~\cite[Theorem~2.1]{zhang2021-1}, we also start from analysing first-order linear differential equation $f'-uf=w$, where $u$ is a nonzero polynomial and $w$ is a meromorphic function with at most finitely many poles. Let $p(z)$ be a primitive function of $u$ and suppose that $\deg(p(z))=k\geq1$. If $f$ is meromorphic, then there is a rational function $v(z)$ such that $v(z)\to 0$ as $z\to\infty$ and $h(z)=f(z)-v(z)$ is entire. It follows that $f(z)=h(z)+v(z)$ and $h$ satisfies $h'-uh=w-(v'-uv)$ and $w-(v'-uv)$ is an entire function. By elementary integration, the meromorphic solutions of $f'-uf=w$ are $f=ce^{p(z)}+H(z)$, where
\begin{equation}\label{Hfunction}
H(z)=e^{p(z)}\int_0^{z}w(t) e^{-p(t)}dt.
\end{equation}
To study the growth behavior of this function, a useful tool is the Phragm\'{e}n--Lindel\"{o}f theorem (see \cite[Theorem~7.3]{hollandasb}): Let $f(z)$ be an analytic function, regular in a region $D$ between two straight lines making an angle $\pi/\tau_1$ at the origin, and on the lines themselves. Suppose that $|f(z)|\leq M$ on the line, and that, as $r\to\infty$ $|f(z)|=O(e^{r^{\tau_2}})$, where $\tau_2<\tau_1$, uniformly in the angle. Then actually $|f(z)|\leq M$ holds throughout the region. Moreover, if $f(z)\to c_1$ and $f(z)\to c_2$ as $z\to\infty$ along the two lines, respectively, then $c_1=c_2$ and $f(z)\to c_1$ uniformly as $z\to\infty$ in $D$. Using the Phragm\'{e}n--Lindel\"{o}f theorem, the present author proved the following

\begin{lemma}[see \cite{zhang2021,zhang2021-1}]\label{first-order lemma}
Let $p(z)$ be a polynomial with degree $k\geq 1$ and $w$ be a nonzero polynomial. Then there is an integer $N$ such that for each $S_j$ where $\delta(p,\theta)>0$, there is a constant $a_j$ such that $|H(re^{i\theta})-a_je^{p(re^{i\theta})}|= O(r^N)$ uniformly in $\overline{S}_{j,\epsilon}$, and for each $S_j$ where $\delta(p,\theta)<0$ and any constant $a$, $|H(re^{i\theta})-ae^{p(re^{i\theta})}|= O(r^N)$ uniformly in $\overline{S}_{j,\epsilon}$.
\end{lemma}

Most arguments we use below are the same as that in the proof of~\cite[Theorem~2.1]{zhang2021-1}. We also first introduce the definition of \emph{$R$--set}: An $R$--set in the complex plane is a countable union of discs whose radii have finite sum. Let $f(z)$ be an entire solution of \eqref{EQ1}. We denote the union of all $R$--sets associated with $f(z)$ and each coefficient of $P(z,f)$ by $\tilde{R}$ from now on. In the proof of Theorem~\ref{maintheorem3}, after taking the derivatives on both sides of equation \eqref{EQ1}, there may be some new coefficients appearing in the resulting equations. We will always assume that $\tilde{R}$ also contains those $R$-sets associated with these new coefficients.

As in the proof of~\cite[Theorem~2.1]{zhang2021-1}, we first reduce \eqref{EQ1} into a non-homogeneous linear differential equation with rational coefficients. Now, with all coefficients of $P(z,f)$ being combinations of functions in $\mathcal{S}$, the key lemma for this aim is the following

\begin{lemma}\label{growthlemmapre}
Under the assumptions of Theorem~\ref{maintheorem3}, $\sigma(f)=k$ and $\alpha$ is real. Moreover, for any $\theta\in[0,2\pi)$ such that the ray $z=re^{i\theta}$ meets finitely discs in $\tilde{R}$,
\begin{itemize}
  \item [(1)]
when $-1\leq \alpha<0$, if $\delta(p_1,\theta)>0$, then $|f(re^{i\theta})^n|=(1+o(1)) |b_1(re^{i\theta})e^{p_1(re^{i\theta})}|$, $r\to\infty$; if $\delta(p_2,\theta)>0$, then $|f(re^{i\theta})^n|=(1+o(1))|b_2(re^{i\theta})e^{p_2(re^{i\theta})}|$, $r\to\infty$;

  \item [(2)]
when $0<\alpha<1$, if $\delta(p_1,\theta)>0$, then $|f(re^{i\theta})^n|=(1+o(1)) |b_1(re^{i\theta})e^{p_1(re^{i\theta})}|$, $r\to\infty$; if $\delta(p_1,\theta)<0$, then there is an integer $N$ such that $|f(re^{i\theta})|\leq r^N$ for all large $r$.
\end{itemize}
\end{lemma}

\renewcommand{\proofname}{Proof of Lemma~\ref{growthlemmapre}.}
\begin{proof}
Since $\alpha\not=1$, then by Steinmetz's result~\cite{Steinmetz1978} for exponential polynomials, we have $T(r,b_1e^{p_1}+b_2e^{p_2})=K(1+o(1))r^k$ for some nonzero constant~$K$ depending only on~$\alpha$.
Recall that the coefficients of equation~\eqref{EQ1} are combinations of functions in $\mathcal{S}$. By the lemma on the logarithmic derivative, we deduce from equation \eqref{EQ1} that
\begin{equation}\label{Hfunctionag;oi}
\begin{split}
T\left(r,b_1e^{p_1}+b_2e^{p_2}\right)&=m\left(r,b_1e^{p_1}+b_2e^{p_2}\right)\\
&=m\left(r,f^n+P(z,f)\right)\leq nm(r,f)+O(\log r).
\end{split}
\end{equation}
Therefore, $f$ is transcendental and $T(r,f)\geq K_1 r^k$ for some positive constant $K_1$. On the other hand, by the lemma on the logarithmic derivative we also have from equation~\eqref{EQ1} that
\begin{equation}\label{Hfunctionag;oi1}
\begin{split}
nT(r,f)&=T\left(r,f^n\right)=m\left(r,f^n\right)=m\left(r,b_1e^{p_1}+b_2e^{p_2}-P(z,f)\right)\\
&\leq m\left(r,b_1e^{p_1}+b_2e^{p_2}\right)+m\left(r,P(z,f)\right)+O(1)\\
&\leq K(1+o(1))r^k+(n-1)m(r,f)+O(\log r),
\end{split}
\end{equation}
which yields that $T(r,f)\leq K_2 r^k$ for some positive constant $K_2$. This together with $T(r,f)\geq K_1 r^k$ yields $\sigma(f)=k$. Then by definition of $\mathcal{S}$ and looking at the proof of \cite[Theorem~2.1]{zhang2021-1}, we see that $\alpha$ is real. Now, $-1\leq \alpha<0$ or $0<\alpha<1$.

Recall that $\theta_1=-\pi/(2k)$ and from \eqref{expon1} that $\delta(p_1,\theta)=\cos k\theta$ and $\delta(p_2,\theta)=\alpha \cos k\theta$. When $\alpha<0$, we see that $\delta(p_1,\theta)$ and $\delta(p_2,\theta)$ have opposite signs for each $\theta$ in the sectors $S_j$ defined in \eqref{expon10} for $p_1$ and $\delta(p_1,\theta)>0$ for $\theta$ in the sectors $S_{j}$ where $j\in J_1$; when $\alpha>0$, we see that $\delta(p_1,\theta)>0$ and $\delta(p_2,\theta)>0$ simultaneously for each $\theta$ in the sectors $S_{j}$ where $j\in J_1$ and $\delta(p_1,\theta)<0$ and $\delta(p_2,\theta)<0$ simultaneously for each $\theta$ in the sectors $S_{j}$ where $j\in J_2$. Then we see that the assertion (1) and the assertion (2) for the case that $\delta(p_1,\theta)>0$ can be obtained by directly following the proof of \cite[Lemma~2.5]{zhang2021}.

Now we consider the growth behavior of $f(z)$ along the ray $z=re^{i\theta}$ such that $\delta(p_1,\theta)<0$ when $0<\alpha<1$. Let $\varepsilon>0$ be given. By~\cite[Corollary~1]{gundersen:88}, there exists a constant $r_0=r_0(\theta)>1$ such that for all $z$ on the ray $z=re^{i\theta}$ which does not meet $\tilde{R}$ when $r\geq r_0$, and for all positive integers $j$,
\begin{equation}\label{lemmagrw0}
\left|\frac{f^{(j)}(re^{i\theta})}{f(re^{i\theta})}\right|\leq r^{j(k-1+\varepsilon)}.
\end{equation}
Since all coefficients of $P(z,f)$ are combinations of functions in $\mathcal{S}$, then for each coefficient of $P(z,f)$, say $a_l$, by~\cite[Corollary~1]{gundersen:88}, we also have, along the ray $z=re^{i\theta}$, that
\begin{equation}\label{lemmagrw1}
\left|a_l(re^{i\theta})\right|\leq r^{M},
\end{equation}
for sufficiently large $r$ and some large integer $M$. Recalling from the introduction that $P(z,f)=\sum^{m}_{l=1}a_{l}f^{n_{l0}}(f')^{n_{l1}}\cdots(f^{(s)})^{n_{ls}}$, where $m$ is an integer and $n_{l0}+n_{l1}+\cdots n_{ls}\leq n-1$, we may write
\begin{equation}\label{lemmaQ10}
P(z,f)=\sum^{m}_{l=1}\hat{a}_{l}f^{n_{l0}+n_{l1}+\cdots+n_{ls}},
\end{equation}
with the new coefficients $\hat{a}_l=a_l(f'/f)^{n_{l1}}\cdots(f^{(s)}/f)^{n_{ls}}$, where $n_{l0},\cdots,n_{ls}$ are nonnegative integers. Note that the greatest order of the derivatives of $f$ in $P(z,f)$ is equal to $s\geq 0$. Suppose now that $|f(r_je^{i\theta})|\geq r_j^{N}$ for some infinite sequence $z_j=r_je^{i\theta}$ and some large $N\geq M+s(k-1+\varepsilon)$. Then, from \eqref{EQ1}, \eqref{lemmagrw0}, \eqref{lemmagrw1} and \eqref{lemmaQ10} we have
\begin{equation}\label{lemmaQ3}
\begin{split}
&\left|b_1(r_je^{i\theta})e^{p_1(r_je^{i\theta})}+b_2(r_je^{i\theta})e^{p_2(r_je^{i\theta})}\right|\\
=\, &\left|f(r_je^{i\theta})^n\right|\left|1+\frac{P(r_je^{i\theta},f(r_je^{i\theta}))}{f(r_je^{i\theta})^n}\right|\geq(1-o(1))r^{nN},
\end{split}
\end{equation}
which is impossible when $r_j$ is large since $b_1(r_je^{i\theta})e^{p_1(r_je^{i\theta})}+b_2(r_je^{i\theta})e^{p_2(r_je^{i\theta})}\to 0$
as $z_j\to\infty$. Therefore, along the ray $z=re^{i\theta}$ such that $\delta(p_1,\theta)<0$ we must have $|f(re^{i\theta})|\leq r^{N}$ for all large $r$ and some integer $N$. Thus our second assertion follows.

\end{proof}

Now we begin to prove Theorem~\ref{maintheorem3}.

\renewcommand{\proofname}{Proof of Theorem~\ref{maintheorem3}.}

\begin{proof}

For simplicity, we denote $P=P(z,f)$. By taking the derivatives on both sides of \eqref{EQ1} and eliminating $e^{p_2}$ and $e^{p_1}$ from \eqref{EQ1} and the resulting equation, respectively, we get the following two equations:
\begin{eqnarray}
b_2B_2f^n-nb_2f^{n-1}f'+b_2B_2 P-b_2P'&=&A_1{e^{p_1}},\label{1Eq2.2}\\
b_1B_1f^n-nb_1f^{n-1}f'+b_1B_1 P-b_1P'&=&-A_1{e^{p_2}},\label{1Eq3.1}
\end{eqnarray}
where $B_1=b_1'/b_1+p_1'$, $B_2=b_2'/b_2+p_2'$ and $A_1=b_1b_2(B_2-B_1)$. Note that $B_1B_2A_1\not\equiv0$. By differentiating on both sides of \eqref{1Eq2.2} and then eliminating $e^{p_1}$ from \eqref{1Eq2.2} and the resulting equation, we get
\begin{equation}\label{1Eq2.3}
h_1 f^n+h_2f^{n-1}f'+h_3f^{n-2}(f')^2+h_4f^{n-1}f''+P_1=0,
\end{equation}
where
$h_1=b_2B_2(A_1'+p_1'A_1)-(b_2B_2)'A_1$, $h_2=-n b_2A_1(p_1'+p_2')-nb_2A_1'$,
$h_3=n(n-1)b_2A_1$, $h_4=nb_2A_1$,
and
$P_1=(A_1'+p_1'A_1)(b_2B_2P-b_2P')-A_1(b_2B_2P-b_2P')'$
is a differential polynomial in $f$ of degree $\leq n-1$. By Lemma~\ref{growthlemmapre} and our assumption, $\alpha$ is a nonzero real number such that $-1\leq \alpha< 1$. Below we consider the two cases where $-1\leq \alpha< 0$ and $0<\alpha<1$, respectively.

\vskip 4pt

\noindent\textbf{Case~1:} $-1\leq \alpha< 0$.

\vskip 4pt

We multiply both sides of equations \eqref{1Eq2.2} and \eqref{1Eq3.1} and obtain
\begin{equation}\label{1Eq3.2}
g_1 f^{2n}+g_2f^{2n-1}f'+g_3f^{2n-2}(f')^2+P_2=-A_1^2e^{p_1+p_2},
\end{equation}
where
$g_1=b_1b_2B_1B_2$, $g_2=-nb_1b_2(B_1+B_2)$, $g_3=n^2b_1b_2$
and
$P_2=b_1b_2(B_2f^n-nf^{n-1}f')(B_1P-P')+b_1b_2(B_1f^n-nf^{n-1}f')(B_2P-P')+b_1b_2(B_1P-P')(B_2P-P')$
is a differential polynomial in $f$ of degree $\leq 2n-1$. By eliminating $(f')^2$ from \eqref{1Eq2.3} and \eqref{1Eq3.2}, we get
\begin{equation}\label{1Eq3.3}
\begin{split}
f^{2n-1}\left[(g_3h_1-h_3g_1)f+(g_3h_2-h_3g_2)f'+g_3h_4f''\right]+P_3=h_3A_1^2e^{p_1+p_2},
\end{split}
\end{equation}
where $P_3=g_3f^nP_1-h_3P_2$ is a differential polynomial in $f$ of degree $\leq 2n-1$. For simplicity, we denote
\begin{equation}\label{1Eq3.9a10a}
\varphi=\frac{h_3A_1^2}{g_3h_4}\frac{e^{p_1+p_2}}{f^{2n-1}}-\frac{1}{g_3h_4}\frac{P_3}{f^{2n-1}}.
\end{equation}
Recalling $B_1=b_1'/b_1+p_1'$ and $B_2=b_2'/b_2+p_2'$, we get from equation \eqref{1Eq3.3} that
\begin{equation}\label{1Eq3.9a1}
f''+H_1f'+H_2f=\varphi,
\end{equation}
where
\begin{equation}\label{1Eq3.9a1fu1}
\begin{split}
H_1&=\frac{h_2}{h_4}-\frac{g_2h_3}{g_3h_4}=-\left[\frac{1}{n}(p_1'+p_2')-\frac{n-1}{n}\left(\frac{b_1'}{b_1}+\frac{b_2'}{b_2}\right)+\frac{A_1'}{A_1}\right],\\
H_2&=\frac{h_1}{h_4}-\frac{g_1h_3}{g_3h_4}=\frac{1}{n}\left[B_2\left(\frac{A_1'}{A_1}-\frac{b_1'}{b_1}\right)-\frac{(b_2B_2)'}{b_2}\right]+\frac{1}{n^2}B_1B_2.
\end{split}
\end{equation}

Now we prove that $\varphi$ is a rational function. Recall that $b_1,b_2,p_1,p_2$ are all polynomials and $B_1=b_1'/b_1+p_1'$, $B_2=b_2'/b_2+p_2'$ and $A_1=b_1b_2(B_2-B_1)$. Since $f$ is entire, we see that $\varphi$ has only finitely many poles. By Lemma~\ref{growthlemmapre}, $\sigma(f)=k$. By the lemma on the logarithmic derivative, we deduce from \eqref{1Eq3.9a1} that
\begin{equation}\label{1Eq3.9a1fu1hyr}
\begin{split}
T(r,\varphi)=m(r,\varphi)+O(\log r)\leq m(r,f)+O(\log r)=T(r,f)+O(\log r).
\end{split}
\end{equation}
Therefore, $\sigma(\varphi)\leq k$. Now let $\theta\in [0,2\pi)$ be such that $\delta(p_1,\theta)\not=0$ and $z=re^{i\theta}$ is a ray that meets only finitely discs in $\tilde{R}$. Since $\alpha<0$, then by Lemma~\ref{growthlemmapre}~(1) we see that in both cases that $\delta(p_1,\theta)>0$ and $\delta(p_1,\theta)<0$ we always have $|e^{p_1(re^{i\theta})+p_2(re^{i\theta})}/f(re^{i\theta})^{2n-1}|\to 0$ as $r\to\infty$ along the ray $z=re^{i\theta}$. Together with~\cite[Corollary~1]{gundersen:88} we see from \eqref{1Eq3.9a10a} that there is some integer $N$ such that $|\varphi(re^{i\theta})|\leq r^{N}$ for all large $r$. Then by the Phragm\'{e}n--Lindel\"{o}f theorem we see that $|\varphi|\leq r^{N}$ uniformly in each $\overline{S}_{j,\epsilon}$, $j=1,2,\cdots,2k$, for some integer $N=N(j)$. Since $\epsilon$ can be arbitrarily small, then by the Phragm\'{e}n--Lindel\"{o}f theorem again we conclude that $\varphi$ is a rational function. From now on we fix one large $N$.

Recall that $B_2=b_2'/b_2+p_2'$. Denote $F_1=f'-(B_1/n)f$. Then by simple computations we obtain from \eqref{1Eq3.9a1} that
\begin{equation}\label{1Eq3.9a10a hfo0}
F_1'-\left(\frac{1}{n}p_2'-\frac{b_1'}{b_1}-\frac{n-1}{n}\frac{b_2'}{b_2}+\frac{A_1'}{A_1}\right)F_1=\varphi.
\end{equation}
Denote $\xi_1=p_2'/n-b_1'/b_1-(n-1)b_2'/nb_2+A_1'/A_1$. Then the general solution of the homogeneous equation $F_1'-\xi_1 F_1=0$ is defined on a finite-sheeted Riemann surface and is of the form $F_1=C_2b_2^{1/n}A_1/(b_1b_2)e^{p_2/n}$, where $C_2$ is a constant and $b_2^{1/n}$ is in general an algebraic function (see \cite{Katajamaki1993algebroid} for the theory of algebroid functions). Suppose that $\Gamma_2$ is a particular solution of $F_1'-\xi_1F_1=\varphi$. We may write the meromorphic solution of this equation as $F_1=C_2b_2^{1/n}A_1/(b_1b_2)e^{p_2/n}+\Gamma_2$. By an elementary series expansion analysis around the zeros of $b_2$, we conclude that $\Gamma_2/b_2^{1/n}$ is a meromorphic function. This implies that $b_2$ is an $n$-square of some polynomial. Then by Lemma~\ref{first-order lemma} we integrate the equation \eqref{1Eq3.9a10a hfo0} along the ray $z=re^{i\theta}$ in $S_2$ such that $\delta(p_2,\theta)>0$ and obtain
\begin{equation}\label{1Eq3.9a10a hfo1}
F_1=f'-\frac{1}{n}B_1f=\frac{c_2}{n}\frac{b_2^{1/n}A_1}{b_1b_2}e^{p_2/n}+\Gamma_2,
\end{equation}
where
\begin{equation}\label{1Eq3.9a10a hfo1  jiq}
\Gamma_2=\frac{A_1b_2^{1/n}}{b_1b_2}e^{p_2/n}\int_{0}^{z} e^{-p_2/n}\frac{b_1b_2}{A_1b_2^{1/n}} \varphi dt-a_{2,2}\frac{A_1b_2^{1/n}}{b_1b_2}e^{p_2/n},
\end{equation}
where $a_{2,2}=a_{2,2}(\theta)$ is a constant such that $|\Gamma_2|=O(r^N)$ along the ray $z=re^{i\theta}$ in $S_2$. Now, for $z\in S_{j,\epsilon}$ where $j\in J_2$, we have $\delta(p_2,\theta)>0$ and so $\Gamma_2=(c_2d_{2,j}/n)b_2^{1/n}A_1/(b_1b_2)e^{p_2/n}+\gamma_{2,j}$, where $d_{2,j}$ are some constants related to a sector $S_{j,\epsilon}$ and $|\gamma_{2,j}|=O(r^N)$ uniformly in $\overline{S}_{j,\epsilon}$. Of course, for $j=2$, we have $d_{2,2}=0$. Furthermore, $|\Gamma_2|=O(r^N)$ uniformly
in $\overline{S}_{j,\epsilon}$ where $j\in J_1$. We then define $d_{2,j}=0$ for $j\in J_1$.

Similarly, denoting that $\xi_2=p_1'/n-b_2'/b_2-(n-1)b_1'/nb_1+A_1'/A_1$ we also have $F_2'-\xi_2F_2=\varphi$ and it follows by integration that $F_2=-(c_1/n)b_1^{1/n}A_1/(b_1b_2)e^{p_1/n}+\Gamma_1$, where
$\Gamma_1=-(c_1d_{1,j}/n)b_1^{1/n}A_1/(b_1b_2)e^{p_1/n}+\gamma_{1,j}$, where $d_{l,j}$ are some constants related to a sector $S_{j,\epsilon}$ and $|\gamma_{1,j}|=O(r^N)$ uniformly in $\overline{S}_{j,\epsilon}$ for $j\in J_1$. Of course, for $j=1$, we have $d_{1,1}=0$. Furthermore, $|\Gamma_1|=O(r^N)$ uniformly in $\overline{S}_{j,\epsilon}$ where $j\in J_2$. We then define $d_{1,j}=0$ for $j\in J_2$.

Denoting $B=n/(B_2-B_1)$, we have $f=B(F_1-F_2)$. Together with the relation $A_1=b_1b_2(B_2-B_1)$, we have $f=c_1b_1^{1/n}e^{p_1/n}+c_2b_2^{1/n}e^{p_2/n}+\eta$ with an entire function $\eta=B(\Gamma_2-\Gamma_1)$. We see that $\eta=c_2d_{2,j}b_2^{1/n}e^{p_2/n}+B(\gamma_{2,j}-\gamma_{1,j})$ when $j\in J_1$ and $\eta=c_1d_{1,j}b_1^{1/n}e^{p_1/n}+B(\gamma_{2,j}-\gamma_{1,j})$ when $j\in J_2$.

Now we determine $d_{1,j}$ and $d_{2,j}$. By~\cite[Corollary~1]{gundersen:88}, we may suppose that along the ray $z=re^{i\theta}$ we have $|f^{(j)}(re^{i\theta})/f(re^{i\theta})|=r^{j(k-1+\varepsilon)}$ for all $j>0$ for all sufficiently large $r$ and thus write $P$ in the form in \eqref{lemmaQ10} with the new coefficients $\hat{a}_l=a_l(f'/f)^{n_{l1}}\cdots(f^{(s)}/f)^{n_{ls}}$, where $n_{l1},\cdots,n_{ls}$ are nonnegative integers. For simplicity, denote $D_{1,j}=c_1+c_1d_{1,j}$. By substituting $f=c_1b_1^{1/n}e^{p_1/n}+c_2b_2^{1/n}e^{p_2/n}+\eta$ into \eqref{EQ1}, we obtain, for $z=re^{i\theta}$ for a $\theta$ in $S_j$ and $j\in J_1$,
\begin{equation}\label{EQ1aa}
\begin{split}
&\left(D_{1,j}^n-1\right)b_1e^{p_1}+\sum_{k_0=1}^{n-1}\binom{n}{k_0}\left(D_{1,j}b_1^{1/n}\right)^{n-k_0}(c_2b_2^{1/n})^{k_0}e^{[(n-k_0)p_1+k_0p_2]/n}\\
&+(c_2^n-1)b_2e^{p_2}+\sum_{s=1}^{n}\sum_{k_s=0}^{n-s}\alpha_{s,k_s}e^{[(n-s-k_s)p_1+k_sp_2]/n}=0,
\end{split}
\end{equation}
where $\alpha_{s,k_s}$, $s=1,\cdots,n$, $k_s=0,\cdots,n-s$, are functions satisfying $|\alpha_{s,k_s}(re^{i\theta})|=O(r^N)$ along the ray $z=re^{i\theta}$. By letting $r\to\infty$ along the above ray $z=re^{i\theta}$ such that $\delta(p_1,\theta)>0$ and comparing the growth on both sides of the above equation we conclude that $c_1^n(1+d_{1,j})^n=1$. Since $d_{1,1}=0$, we have $c_1^n=1$ and $d_{1,j}=\mu_{1,j}-1$ for some $\mu_{1,j}$ such that $\mu_{1,j}^n=1$. Similarly, we can prove that $d_{2,j}=\mu_{2,j}-1$ for some $\mu_{2,j}$ such that $\mu_{2,j}^n=1$. In particular, when $k=1$, since $d_{1,1}=d_{2,2}=0$ and $|\eta_j|=O(r^N)$ uniformly in the sectors $\overline{S}_{j,\epsilon}$, $j=1,2$ and since $\epsilon$ can be arbitrarily small, by the Phragm\'{e}n--Lindel\"{o}f theorem we conclude that $\eta$ is a polynomial. Thus our first assertion follows.

\vskip 4pt

\noindent\textbf{Case~2:} $0<\alpha<1$.

\vskip 4pt

As in the proof of \cite[Theorem~2.1]{zhang2021-1}, we first define some functions in the following way: We let $m$ be the smallest integer such that $\alpha \leq[(m+1)n-1]/[(m+1)n]$ and ${\iota_0,\cdots,\iota_m}$ be a finite sequence of functions such that
\begin{equation}\label{1Eq3.27a prepare0}
\begin{split}
\iota_0&=\frac{A_1}{nb_1},\\
\iota_j&=(-1)^{j}\left(\frac{A_1}{nb_1}\right)^{j+1}(jn-1)\cdots(n-1), \quad j=1,2,\cdots,m.
\end{split}
\end{equation}
Recall that $B_1=b_1'/b_1+p_1'$. We also let ${\kappa_0,\cdots,\kappa_m}$ be a finite sequence of functions defined in the following way:
\begin{equation}\label{1Eq3.27a prepare}
\begin{split}
\kappa_0&=\frac{1}{n}\frac{b_1'}{b_1}+\frac{1}{n}p_1',\\
\kappa_j&=\frac{\iota_{j-1}'}{\iota_{j-1}}-\frac{jn-1}{n}\frac{b_1'}{b_1}+\left[j(\alpha-1)+\frac{1}{n}\right]p_1', \quad j=1,2,\cdots,m.
\end{split}
\end{equation}
Then we define $m+1$ functions $G_0$, $G_1$, $\cdots$, $G_m$ in the way that $G_0=f'-\kappa_0f$, $G_{1}=G_{0}'-\kappa_1G_{0}$, $\cdots$, $G_{m}=G_{m-1}'-\kappa_mG_{m-1}$. Now we have equation \eqref{1Eq3.1} and it follows that
\begin{equation}\label{1Eq3.27a}
G_0=f'-\kappa_0f=\iota_0\frac{{e^{p_2}}}{f^{n-1}}+W_0,
\end{equation}
where $W_0=-(B_1P-P')/(nf^{n-1})$. Moreover, when $m\geq 1$, by simple computations we obtain
\begin{eqnarray*}
G_1&=&G_0'-\kappa_1G_0=\iota_1\frac{e^{2p_2}}{f^{2n-1}}+W_1,\\
W_1&=&W'_0-\kappa_1W_0-(n-1)\iota_0\frac{e^{p_2}}{f^n}W_0,
\end{eqnarray*}
and by induction we obtain
\begin{eqnarray}
G_{j}&=&G_{j-1}'-\kappa_jG_{j-1}=\iota_j\frac{e^{(j+1)p_2}}{f^{(j+1)n-1}}+W_{j}, \quad j=1,\cdots,m,\label{1Eq3.27}\\
W_{j}&=&W'_{j-1}-\kappa_jW_{j-1}-(jn-1)\iota_{j-1}\frac{e^{jp_2}}{f^{jn}}W_0, \quad j=1,\cdots,m.\label{1Eq3.27a13}
\end{eqnarray}
For an integer $l\geq0$, by elementary computations it is easy to show that $W_0^{(l)}=W_{0l}/f^{n+l-1}$, where $W_{0l}=W_{0l}(z,f)$ is a differential polynomial in $f$ of degree $\leq n+l-1$, and also that $(e^{p_2}/f^{n})^{(l)}=e^{p_2}W_{1l}/f^{n+l}$, where $W_{1l}=W_{1l}(z,f)$ is a differential polynomial in $f$ of degree $\leq n+l$. We see that $W_j$, $1\leq j\leq m$, is formulated in terms of $W_0$ and $e^{p_2}/f^{n}$ and their derivatives. We may write
\begin{equation}\label{1Eq3.28a1}
G_m=\iota_m\frac{e^{(m+1)p_2}}{f^{(m+1)n-1}}+F(W_0,e^{p_2}/f^{n}),
\end{equation}
where $F(W_0,e^{p_2}/f^{n})$ is a combination of $W_0$ and $e^{p_2}/f^{n}$ and their derivatives with functions being combinations of functions in $\mathcal{S}$. Moreover, from the recursion formula $G_j=G_{j-1}'-\kappa_jG_{j-1}$, $j\geq 1$, and $G_0=f'-\kappa_0f$, we easily deduce that $f$ satisfies the linear differential equation
\begin{equation}\label{1Eq3.28lin}
f^{(m+1)}-\hat{t}_{m}f^{(m)}+\cdots+(-1)^{m+1}\hat{t}_0f=G_m,
\end{equation}
where $\hat{t}_m$, $\hat{t}_{m-1}$, $\cdots$, $\hat{t}_0$ are functions formulated in terms of $\kappa_0$, $\cdots$, $\kappa_m$ and their derivatives.

Now we prove that $G_m$ is a rational function. Recall that $b_1,b_2,p_1,p_2$ are all polynomials. Since $f$ is entire, then by the definitions of $\kappa_0$ and $\kappa_j$ in \eqref{1Eq3.27a prepare}, we see that $G_m$ has only finitely many poles. With an application of the lemma on the logarithmic derivative as in previous case, we deduce from \eqref{1Eq3.28lin} that $\sigma(G_m)\leq \sigma(f)=k$. Now let $\theta\in [0,2\pi)$ be such that $\delta(p_1,\theta)\not=0$ and $z=re^{i\theta}$ be a ray that meets only finitely may discs in $\tilde{R}$. By~\cite[Corollary~1]{gundersen:88} and Lemma~\ref{growthlemmapre}~(2), we see from \eqref{1Eq3.28lin} that there is some integer $N$ such that $|G_m(re^{i\theta})|\leq r^{N}$ for all large $r$ along the ray $z=re^{i\theta}$ such that $\delta(p_1,\theta)<0$. On the other hand, by Lemma~\ref{growthlemmapre}~(2) there is some integer $N$ such that
\begin{itemize}
  \item [(1)]
if $\alpha<[(m+1)n-1]/[(m+1)n]$, then $|e^{(m+1)p_2(re^{i\theta})}/f(re^{i\theta})^{(m+1)n-1}|\to 0$ as $r\to\infty$ along the ray $z=re^{i\theta}$ such that $\delta(p_1,\theta)>0$;

  \item [(2)]
if $\alpha=[(m+1)n-1]/[(m+1)n]$, then $|e^{(m+1)p_2(re^{i\theta})}/f(re^{i\theta})^{(m+1)n-1}|\leq e^{Nr^{k-1}}$ for all large $r$ along the ray $z=re^{i\theta}$ such that $\delta(p_1,\theta)>0$.

\end{itemize}
Note that $e^{p_2(re^{i\theta})}/f(re^{i\theta})^{n}\to 0$ as $r\to \infty$ along the ray $z=re^{i\theta}$ such that $\delta(p_1,\theta)>0$. In case~(1), together with~\cite[Corollary~1]{gundersen:88} we see from \eqref{1Eq3.28a1} that $|G_m(re^{i\theta})|\leq r^{N}$ for all large $r$ and thus by the Phragm\'{e}n--Lindel\"{o}f theorem we see that $|G_m|\leq r^{N}$ uniformly in each $\overline{S}_{j,\epsilon}$, $j\in J_2$, for some integer $N=N(j)$; in case~(2), together with~\cite[Corollary~1]{gundersen:88} we see from \eqref{1Eq3.28a1} that $|G_m(re^{i\theta})|\leq e^{Nr^{k-1}}$ for all large $r$ and, since the set of rays $z=re^{i\theta}$ meeting infinitely many discs in $\tilde{R}$ has zero linear measure, then by the Phragm\'{e}n--Lindel\"{o}f theorem we see that $|G_m|\leq e^{Nr^{k-1}}$ uniformly in each $\overline{S}_{j,\epsilon}$, $j\in J_2$, for some integer $N=N(j)$. Since $\epsilon$ can be arbitrarily small, then in either case of (1) and (2) by the Phragm\'{e}n--Lindel\"{o}f theorem again we conclude that $G_m$ is a rational function. From now on we fix one large~$N$.

We denote $D_0=b_1^{1/n}$ and $D_j=\iota_{j-1}b_1^{-j}b_1^{1/n}$, $j=1,\cdots,m$. Now we choose one $\theta$ such that $\delta(p_1,\theta)>0$ and let $z=re^{i\theta}\in S_1$. Let $t_0=1/n$, $t_1=(\alpha-1)+1/n$, $\cdots$, $t_m=m(\alpha-1)+1/n$. Similarly as in the proof of \cite[Theorem~2.1]{zhang2021-1}, we may use Lemma~\ref{first-order lemma} to integrate the recursion formulas $G_j=G_{j-1}'-\kappa_jG_{j-1}$ from $j=m$ to $j=1$ along the above ray $z=re^{i\theta}$ such that $\delta(p_1,\theta)>0$ inductively and finally integrating $G_0=f'-\kappa_0f$ along this ray $z=re^{i\theta}$ to obtain
\begin{equation}\label{1Eq3.28a2}
f=b_1^{1/n}\sum_{i=0}^mc_i\left(\frac{b_2}{b_1}\right)^ie^{t_ip_1}+H_0,
\end{equation}
where $c_0$, $\cdots$, $c_m$ are constants and
\begin{equation}\label{1Eq3.28a5}
\begin{split}
H_{0}=b_1^{1/n}e^{t_0p_1}\int_{0}^{z}b_1^{-1/n}e^{-t_0p_1}H_1ds-a_{0}b_1^{1/n}e^{t_0p_1},
\end{split}
\end{equation}
where $a_0=a_{0}(\theta)$ is a constant such that $|H_0|=O(r^N)$ along the ray $z=re^{i\theta}$.

As is shown in the proof of \cite[Theorem~2.1]{zhang2021-1}, $b_1$ is an $n$-square of some polynomial and we may write the entire solution of \eqref{EQ1} as $f=\gamma_1\sum_{j=0}^mc_j(b_2/b_1)^je^{t_jp_1}+\eta$, where $\gamma_1$ is a polynomial such that $\gamma_1^n=b_1$ and $\eta$ is a meromorphic function with at most finitely many poles. Then we can integrate $G_j=G_{j-1}'-\kappa_jG_{j-1}$ from $j=m$ to $j=1$ inductively and finally integrate $G_0=f'-\kappa_0f$ to obtain that $H_0$ is a meromorphic function with at most finitely many poles. We choose $\eta=H_0$. Recall that along the ray $z=re^{i\theta}$ such that $\delta(p_1,\theta)>0$ and $z=re^{i\theta}\in S_1$, we have $|H_0|=O(r^N)$. Denote $g=b_1^{1/n}\sum_{i=0}^mc_i(b_2/b_1)^ie^{t_ip_1}$. Then
\begin{equation}\label{1Eq3.9a1hoanr}
\begin{split}
g^n=b_1\sum_{k_0=0}^{mn}C_{k_0}\left(\frac{b_2}{b_1}\right)^{k_0}e^{(k_0t-k_0+1)p_1},
\end{split}
\end{equation}
where
\begin{equation}\label{binocoe}
\begin{split}
C_{k_0}=\sum_{\substack{j_0+\cdots+j_m=n,\\j_1+\cdots+mj_m=k_0}}\frac{n!}{j_0!j_1!\cdots j_m!}c_0^{j_0}c_1^{j_1}\cdots c_m^{j_m}, \quad k_0=0,1,\cdots,mn.
\end{split}
\end{equation}
By~\cite[Corollary~1]{gundersen:88}, we may suppose that along the ray $z=re^{i\theta}$ we have $|f(re^{i\theta})^{(j)}/f(re^{i\theta})|=r^{j(k-1+\varepsilon)}$ for all $j>0$ and all sufficiently large $r$.
By writing $P$ in the form in \eqref{lemmaQ10} with the new coefficients $\hat{a}_l=a_l(f'/f)^{n_{l1}}\cdots(f^{(s)}/f)^{n_{ls}}$, where $n_{l1},\cdots,n_{ls}$ are nonnegative integers, and using~\cite[Corollary~1]{gundersen:88}, we see that each term in $P(z,f)$ of degree $n-j$, $1\leq j\leq n-1$, equals a linear combination of exponential functions of the form $e^{[nk_{j}(\alpha -1)+n-j]p_1/n}$, $0\leq k_{j}\leq (n-j)m$, with coefficients $\beta_j$ having polynomial growth along the ray $z=re^{i\theta}$. Therefore, by substituting $f=g+H_0$ into \eqref{EQ1} we obtain by the same arguments in the proof of \cite[Theorem~2.1]{zhang2021-1} that $c_0^n=1$ when $m=0$, and $c_0^{n}=1$, $nc_0^{n-1}c_1=1$ and $p_2=\alpha p_1$ when $m=1$ and further that $C_{k_0}\equiv0$ for all $2\leq k_0\leq m$ when $m\geq 2$.

Now, $b_1^{1/n}$ denotes a polynomial. By the definition of $\iota_j$ and $D_j$, we see that $D_j$ are rational functions. Recall that $G_m$ is a rational function. By Lemma~\ref{first-order lemma} and looking at the calculations to obtain $H_0$ in \eqref{1Eq3.28a5}, we have, for $z\in S_{j,\epsilon}$, $j\in J_1$, such that $\delta(p_1,\theta)>0$, $H_0=\gamma_1\sum_{l=0}^md_{l,j}(b_2/b_1)^je^{t_jp_1}+\eta_{j}$, where $d_{l,j}$, $l=0,\cdots,m$, are some constants related to a sector $S_{j,\epsilon}$ and $|\eta_j|=O(r^N)$ uniformly in $\overline{S}_{j,\epsilon}$, $j\in J_1$. Of course, for $j=1$, we have $d_{l,1}=0$ for all $l$. Since $c_0^{n}=1$ when $m=0$, $c_0^{n}=nc_0^{n-1}c_1=1$ when $m=1$, and $c_0^{n}=nc_0^{n-1}c_1=1$ and $C_{k_0}=0$ for all $2\leq k_0\leq m$ when $m\geq 2$, then by simple computations, we deduce that $c_j=s_jc_0$ for some nonzero rational numbers $s_j$, $j=0,1,\cdots,m$. Therefore, by considering the growth of $f$ along the ray $z=re^{i\theta}$ such that $z\in S_{j,\epsilon}$, $j\in J_1$, as for the ray $z=re^{i\theta}\in S_1$, we have $(c_0+d_{0,j})^n=1$ when $m=0$, $(c_0+d_{0,j})^n=n(c_0+d_{0,j})^{n-1}(c_1+d_{1,j})=1$ when $m=1$ and further that $\hat{C}_{k_0}=\sum_{\substack{j_0+\cdots+j_m=n,\\j_1+\cdots+mj_m=k_0}}\frac{n!}{j_0!j_1!\cdots j_m!}(c_0+d_{0,j})^{j_0}(c_1+d_{1,j})^{j_1}\cdots (c_m+d_{m,j})^{j_m}=0$ for $k_0=2,\cdots,m$ when $m\geq 2$.
Therefore, for each $j\in J_1$, there is a $\mu_j$ satisfying $\mu_j^n=1$ such that $c_l+a_{l,j}=\mu_j c_l$ for all $l$. Note that $\mu_{1}=1$. Also, we have $|\eta|=O(r^N)$ uniformly in the sectors $\overline{S}_{j,\epsilon}$, $j\in J_2$. In conclusion, we may write $\eta=\gamma_1\sum_{l=0}^m(\mu_j-1)c_{j}(b_2/b_1)^je^{[jn(\alpha-1)+1]p_1/n}+\eta_{j}$, where $\mu_j$ are the $n$-th roots of~$1$ such that $\mu_j=1$, $j=\{1\}\cup\in J_2$, and $|\eta_j|=O(r^N)$ uniformly in the sector $\overline{S}_{j,\epsilon}$. In particular, when $k=1$, since $\epsilon$ can be arbitrarily small, then by the Phragm\'{e}n--Lindel\"{o}f theorem we conclude that $\eta$ is a rational function. This completes the proof.

\end{proof}

\section{An oscillation question of Ishizaki}\label{se3: an oscillation question of Ishizaki}

Let $b_1(z)$, $b_2(z)$ and $b_3(z)$ be three polynomials such that $b_1b_2\not\equiv0$ and $p_1(z)$ and $p_2(z)$ be two polynomials of the same degree~$k\geq 1$ with distinct leading coefficients $1$ and $\alpha$, respectively, and $p_1(0)=p_2(0)=0$. In this section, we use Theorem~\ref{maintheorem3} to investigate the oscillation of the second-order linear differential equation:
\begin{equation}\label{Bank-laine se3}
f''-\left[b_1(z)e^{p_1(z)}+b_2(z)e^{p_2(z)}+b_3(z)\right]f=0.
\end{equation}
There have been several results about the oscillation of equation \eqref{Bank-laine se3} and recently second-order linear differential equations with exponential polynomials are taken into more consideration in \cite{Heittokangasilt2019,Heittokangasilt2021}. The results of Bank, Laine and langely~\cite{banklainelangely1989}, Ishizaki and Kazuya~\cite{Ishizakikazuya19971} and Ishizaki~\cite{Ishizaki19970} can be summarized as follows:
\begin{itemize}
\item [(1)]
if $\alpha$ is non-real, then all nontrivial solutions of \eqref{Bank-laine se3} satisfy $\lambda(f)=\infty$;

\item [(2)]
if $\alpha$ is negative, then all nontrivial solutions of \eqref{Bank-laine se3} satisfy $\lambda(f)=\infty$;

\item [(3)]
if $0<\alpha<1/2$ or if $b_3\equiv0$ and $3/4<\alpha<1$, then all nontrivial solutions of \eqref{Bank-laine se3} satisfy $\lambda(f)\geq k$.
\end{itemize}
Theorem~\ref{maintheorem1} shows that the condition $b_3\equiv0$ in the third result can be removed. Ishizaki~\cite{Ishizaki19970} asked if the third result $\lambda(f)\geq k$ above can be replaced by $\lambda(f)=\infty$. With Theorem~\ref{maintheorem3} at our disposal, we are able to answer this question partially. We prove the following

\begin{theorem}\label{maintheorem4}
Let $0<\alpha<1$ and $m$ be the smallest integer such that $\alpha\leq[2(m+1)-1]/[2(m+1)]$. Suppose that $b_3\equiv0$ in \eqref{Bank-laine se3}. If \eqref{Bank-laine se3} admits a nontrivial solution~$f$ such that $\lambda(f)<\infty$, then $\alpha=[2(m+1)-1]/[2(m+1)]$ and $p_2=\alpha p_1$.
\end{theorem}

We will mainly use the techniques in \cite{Banklangley1987} (see also \cite[Theorem~5.7]{Laine1993}) to prove Theorem~\ref{maintheorem4}. Since $\alpha$ is a positive number, we have $\int_{1}^{\infty}r|A(re^{i\theta})|dr<\infty$ along the ray $z=re^{i\theta}$ such that $\delta(p_1,\theta)<0$.
The following lemma can be proved similarly as in \cite[Lemma~5.16]{Laine1993} by using Gronwall's lemma (see \cite[p.~86]{Laine1993}).

\begin{lemma}\label{lemma4 prepare0}
Under the assumptions of Theorem~\ref{maintheorem4}, all solutions of equation \eqref{Bank-laine se3} satisfy $|f(re^{i\theta})|=O(r)$ as $r\to\infty$ along the ray $z=re^{i\theta}$ such that $\delta(p_1,\theta)<0$.
\end{lemma}

Now we begin to prove Theorem~\ref{maintheorem4}.

\renewcommand{\proofname}{Proof of Theorem~\ref{maintheorem4}.}
\begin{proof}
Let $f$ be a nontrivial solution of equation \eqref{Bank-laine se3} such that $\lambda(f)<\infty$. By Hadamard's factorization theorem we may write $f=\kappa e^{h}$, where $h$ is an entire function and $\kappa$ is the canonical product from the zeros of $f$ satisfying $\rho(\kappa)=\lambda(\kappa)<\infty$. Denoting $g=h'$, then from \eqref{Bank-laine se3} we have
\begin{equation}\label{second tumura-clunie}
g^2+g'+2\frac{\kappa'}{\kappa}g+\frac{\kappa''}{\kappa}=b_1(z)e^{p_1(z)}+b_2(z)e^{p_2(z)}.
\end{equation}
Below we consider the two cases where $0<\alpha\leq 1/2$ and $(2m-1)/(2m)< \alpha\leq [2(m+1)-1]/[2(m+1)]$, $m\geq 1$, respectively.

\vskip 4pt

\noindent{\bf Case 1:} $0<\alpha\leq 1/2$.

\vskip 4pt

By Theorem~\ref{maintheorem3}, we may write $g=\gamma_1e^{p_1/2}+\eta$, where $\gamma_1$ is a polynomial such that $\gamma_1^2=b_1$ and $\eta$ is an entire function such that $|\eta|=O(r^N)$ uniformly in $\overline{S}_{1,\epsilon}$ and $\overline{S}_{2,\epsilon}$. By substituting this expression into equation \eqref{second tumura-clunie}, we obtain
\begin{equation}\label{1Eq3.29 prev0}
\begin{split}
2\gamma_1\left(\frac{\kappa'}{\kappa}+\frac{1}{2}\frac{\gamma_1'}{\gamma_1}+\frac{p_1'}{4}+\eta\right)e^{p_1/2}-b_2e^{p_2}+\frac{\kappa''}{\kappa}+2\eta\frac{\kappa'}{\kappa}+\eta^2+\eta'=0.
\end{split}
\end{equation}
Suppose that $0<\alpha<1/2$. We define
\begin{equation}\label{1Eq3.29 prev0dss}
\begin{split}
w=\kappa \gamma_1^{1/2}e^{p_1/4+\int_{z_0}^{z} \eta dt},
\end{split}
\end{equation}
where $z_0$ is chosen so that $|z_0|$ is large. Then $w$ is analytic outside a finite disc centered at $0$ and satisfy
\begin{equation}\label{1Eq3.29 prev0dss;iajgag}
\begin{split}
\frac{w'}{w}=\frac{\kappa'}{\kappa}+\frac{1}{2}\frac{\gamma_1'}{\gamma_1}+\frac{p_1'}{4}+\eta.
\end{split}
\end{equation}
Dividing by $2\gamma_1e^{p_1/2}$ on both sides of equation \eqref{1Eq3.29 prev0} and then considering the growth of $w'/w$ along the ray $z=re^{i(\theta_2-\epsilon)}$ such that $w$ has no zero around the neighborhood of the ray $z=re^{i(\theta_2-\epsilon)}$, we have by \cite[Corollary~1]{gundersen:88} that $|w'(re^{i\theta})/w(re^{i\theta})|=O(r^{-2})$ as $r\to\infty$. By integration, we obtain that $w(re^{i(\theta_2-\epsilon)})\to a$ as $r\to\infty$ along the ray $z=re^{i(\theta_2-\epsilon)}$ for some nonzero constant $a=a(\theta_2,\epsilon)$. On the other hand, by applying Lemma~\ref{lemma4 prepare0} to equation \eqref{Bank-laine se3} we have $|f(re^{i(\theta_2+\epsilon)})|=O(r)$ along the ray $z=re^{i(\theta_2+\epsilon)}$. Recalling that $f=\kappa e^{h}$ and $g=h'=\gamma_1e^{p_1/2}+\eta$, we may write
\begin{equation}\label{1Eq3.29 prev0dss1pp}
\begin{split}
w=fe^{-h}\gamma_1^{1/2}e^{p_1/4+\int \eta dz}=f\gamma_1^{1/2}e^{p_1/4-\int_{z_0}^z \gamma_1e^{p_1/2}dt}.
\end{split}
\end{equation}
Since $\delta(p_1,\theta_2+\epsilon)<0$ and thus along the ray $z=re^{i(\theta_2+\epsilon)}$ we have $\int_{z_0}^z \gamma_1e^{p_1/2}dt\to c$ for some constant $c=c(\theta_2,\epsilon)$, we see from \eqref{1Eq3.29 prev0dss1pp} that $w$ defined in \eqref{1Eq3.29 prev0dss} satisfies $w(re^{i(\theta+\epsilon)})\to 0$ as $r\to\infty$. Denote
\begin{equation}\label{1Eq3.29 prev0dss1ppoir}
\begin{split}
S_{\epsilon}=\{re^{i\theta}:\theta_2-\epsilon\leq \theta\leq \theta_2+\epsilon\}.
\end{split}
\end{equation}
By choosing $\epsilon$ to be small and applying the Phragm\'{e}n--Lindel\"{o}f theorem to $w$ defined in \eqref{1Eq3.29 prev0dss} in the sector in \eqref{1Eq3.29 prev0dss1ppoir}, we get $a=0$, a contradiction. Therefore, we must have $\alpha=1/2$ when $b_3\equiv0$.

Now, if $k=1$, then obviously $p_2=p_1/2$ since we have assumed $p_1(0)=p_2(0)=0$. If $k>1$, then by Theorem~\ref{maintheorem3} we have $g=\mu_j\gamma_1e^{p_1/2}+\eta_j$, where $\mu_j^n=1$, $\gamma_1$ is a polynomial such that $\gamma_1^2=b_1$ and $\eta_j$ is an entire function such that $|\eta_j|=O(r^N)$ uniformly in $\overline{S}_{j,\epsilon}$. Note that $\eta_j$ has finite order. Denoting $p_3=p_2-p_1/2$, we rewrite equation \eqref{1Eq3.29 prev0} as
\begin{equation}\label{1Eq3.29 prev1}
\begin{split}
\left[b_2e^{p_3}-2\mu_j\gamma_1\left(\frac{\kappa'}{\kappa}+\frac{1}{2}\frac{\gamma_1'}{\gamma_1}+\frac{p_1'}{4}+\eta_j\right)\right]e^{p_1/2}=\frac{\kappa''}{\kappa}+2\eta_j\frac{\kappa'}{\kappa}+\eta_j^2+\eta_j'.
\end{split}
\end{equation}
If $p_2\not\equiv p_1/2$, then $p_3$ is a nonconstant polynomial such that $\deg(p_3)\leq \deg(p_2)-1$. By the definition of $S_j$ in \eqref{expon10}, we may choose a $\theta\in[0,2\pi)$ so that the ray $z=re^{i\theta}$ meets only finitely discs in $\tilde{R}$ and also that $\log |e^{p_1/2}|$ and $\log |e^{p_2-p_1/2}|$ both increase along the ray $z=re^{i\theta}$. By~\cite[Corollary~1]{gundersen:88} we see that $\kappa'/\kappa+\gamma_1'/2\gamma_1+p_1'/4+\eta_j$ and $\kappa''/\kappa+2\eta_j\kappa'/\kappa+\eta_j^2+\eta_j'-b_3$ both have polynomial growth along the ray $z=re^{i\theta}$. Then by comparing the growth on both sides of equation \eqref{1Eq3.29 prev1} along the ray $z=re^{i\theta}$, we get a contradiction. Therefore, we must have $p_2\equiv p_1/2$ when $\alpha=1/2$.

\vskip 4pt

\noindent{\bf Case 2:} $(2m-1)/(2m)<\alpha\leq [2(m+1)-1]/[2(m+1)]$, $m\geq1$.

\vskip 4pt

In this case, by Theorem~\ref{maintheorem3} we already have $p_2=\alpha p_1$ and we may write $g=\gamma_1\sum_{j=0}^mc_j(b_2/b_1)^je^{[2j(\alpha-1)+1]p_1/2}+\eta$, where $m\geq 1$, $\gamma_1$ is a polynomial such that $\gamma_1^2=b_1$ and $\eta$ is a mermorphic function with at most finitely many poles such that $|\eta|=O(r^N)$ uniformly in $\overline{S}_{1,\epsilon}$ and $\overline{S}_{2,\epsilon}$. By substituting this expression into equation \eqref{second tumura-clunie}, we obtain
\begin{equation}\label{1Eq3.29 prev}
\begin{split}
&2\gamma_1\sum_{j=0}^mc_j\left(\frac{b_2}{b_1}\right)^{j}\left[\frac{\kappa'}{\kappa}+\frac{1}{2}\frac{\gamma_1'}{\gamma_1}+j\frac{(b_2/b_1)'}{b_2/b_1}+\frac{2j(\alpha-1)+1}{4}p_1'+\eta\right]e^{L_jp_1}\\
&\gamma_1^2\sum_{k_0=m+1}^{2m}C_{k_0}\left(\frac{b_2}{b_1}\right)^{k_0}e^{M_{k_0}p_1}+\frac{\kappa''}{\kappa}+2\eta\frac{\kappa'}{\kappa}+\eta^2+\eta'=0,
\end{split}
\end{equation}
where $L_j=[2j(\alpha-1)+1]/2$, $M_{k_0}=k_0\alpha-k_0+1$ and the coefficients $C_{k_0}=\sum_{\substack{j_0+\cdots+j_m=2,\\j_1+\cdots+mj_m=k_0}}\frac{2!}{j_0!j_1!\cdots j_m!}c_0^{j_0}c_1^{j_1}\cdots c_m^{j_m}$, $k_0=m+1,\cdots,2m$. Suppose that $\alpha<[2(m+1)-1]/2(m+1)$.
Then, for $k_0=m+1+j$, $j=0,1,\cdots,m-1$, we have $L_{j+1}<k_0\alpha-\alpha+1<L_j$. As in previous case, we also define the function $w$ in \eqref{1Eq3.29 prev0dss},
where $z_0$ is chosen so that $|z_0|$ is large and $w$ is analytic outside a finite disc centered at $0$. It follows that $w'/w$ has the form in \eqref{1Eq3.29 prev0dss;iajgag}. Similarly as in previous case, we first divide by $2c_0\gamma_1e^{p_1/2}$ on both sides of equation \eqref{1Eq3.29 prev} and conclude that $w(re^{i(\theta_2-\epsilon)})\to a$ as $r\to\infty$ along the ray $z=re^{i(\theta_2-\epsilon)}$ for some nonzero constant $a=a(\theta_2,\epsilon)$; then we use the expression $g=h'=\gamma_1\sum_{j=0}^mc_j(b_2/b_1)^je^{[2j(\alpha-1)+1]p_1/2}+\eta$ to derive from \eqref{1Eq3.29 prev0dss} that $w(re^{i(\theta+\epsilon)})\to 0$ as $r\to\infty$ along the ray $z=re^{i(\theta_2+\epsilon)}$. An application of the Phragm\'{e}n--Lindel\"{o}f theorem to $w$ in the sector in \eqref{1Eq3.29 prev0dss1ppoir} then yields a contradiction. We omit those details. Therefore, we must have $\alpha=[2(m+1)-1]/2(m+1)$. We complete the proof.

\end{proof}

\section{Equation~\eqref{bank-laine0} with periodic coefficients in \eqref{bank-laine012}}\label{se4: equation with periodic coefficients}

As mentioned in the introduction, all nontrivial solutions of the second-order linear differential equation $f''+(e^z-b)f=0$ such that $\lambda(f)<\infty$ are given explicit expressions. In this section we solve nontrivial solutions such that $\lambda(f)<\infty$ of the second-order linear differential equation:
\begin{equation}\label{Bank-laine se3simple}
f''-\left(e^{lz}+b_2e^{sz}+b_3\right)f=0,
\end{equation}
where $l$ and $s$ are relatively prime integers such that $l>s\geq 1$, $b_2$ and $b_3$ are constants and $b_2\not=0$. We remark that by using the method in~\cite{langely1986}, we may prove that all nontrivial solutions of equation \eqref{Bank-laine se3simple} satisfy $\lambda(f)=\infty$ when $b_3$ is replaced by a nonconstant polynomial.

Suppose that equation \eqref{Bank-laine se3simple} admits a nontrivial solution such that $\lambda(f)<\infty$. Then $f$ has the form in \eqref{bank-laine01hf1} or \eqref{bank-laine01hf2}. Also, we may write $f=\kappa e^{h}$, where $h$ is an entire function and $\kappa$ is the canonical product from the zeros of $f$ satisfying $\sigma(\kappa)=\lambda(\kappa)<\infty$. Thus we may suppose that $\kappa$ equals a polynomial in $e^{z/2}$ or $e^{z}$ and $h'$ equals a polynomial in $e^{z/2}$ or $e^{z}$. By denoting $g=h'$, from \eqref{Bank-laine se3simple} we have
\begin{equation}\label{second tumura-clunieja;gij0}
g^2+g'+2\frac{\kappa'}{\kappa}g+\frac{\kappa''}{\kappa}=e^{lz}+b_2e^{sz}+b_3.
\end{equation}
By Theorem~\ref{maintheorem3}, we may determine the coefficients $c_j$ in \eqref{bank-laine01hf1} or \eqref{bank-laine01hf2} from equation \eqref{second tumura-clunieja;gij0}. Our main result is the following

\begin{theorem}\label{maintheorem5}
Let $b_2$ and $b_3$ be constants such that $b_2\not=0$ and $l$, $s$ be relatively prime integers such that $l>s\geq1$. Suppose that~\eqref{Bank-laine se3simple} admits two linearly independent solutions $f_1$ and $f_2$ such that $\max\{\lambda(f_1),\lambda(f_2)\}<\infty$. Then $s=1$ and $l=2$.
\end{theorem}

Recall the following well-known result due to Wittich \cite{Wittich1967}. We say that a function $f$ is \emph{subnormal} if $\limsup_{r\to\infty} \log T(r,f)/r=0$. This lemma gives the form of subnormal solutions of second-order linear differential equations with certain periodic functions as coefficients.

\begin{lemma}\label{lemma  solution0pre}
Let $P(z)$ and $Q(z)$ be polynomials in $z$ and not both constants. If $w\not\equiv0$ is a subnormal solution of equation
\begin{equation}\label{wittichequation}
w''+P(e^z)w'+Q(e^z)w=0,
\end{equation}
then $w$ must have the form $w=e^{cz}\left(a_0+a_1e^{z}+\cdots+a_ke^{kz}\right)$, where $k\geq 0$ is an integer and $c$, $a_0$, $\cdots$, $a_k$ are constants with $a_0\not=0$ and $a_k\not=0$. Moreover, we have $c^2+cP(0)+Q(0)=0$.
\end{lemma}

\renewcommand{\proofname}{Proof of Lemma~\ref{lemma  solution0pre}.}
\begin{proof}
By Wittich \cite{Wittich1967}, we have $w=e^{cz}(a_0+a_1e^{z}+\cdots+a_ke^{kz})$. By taking the derivatives of $w$ and then dividing $w'$ and $w''$ by $w$, respectively, we get
\begin{eqnarray}
\frac{w'}{w}&=&\frac{\sum_{j=0}^k(c+j)a_je^{jz}}{\sum_{j=0}^ka_je^{jz}},\label{second tumura-cluniesimple4 fube}\\
\frac{w''}{w}&=&\frac{\sum_{j=0}^k(c^2+2jc+j^2)a_je^{jz}}{\sum_{j=0}^ka_je^{jz}}.\label{second tumura-cluniesimple4 fu1be1}
\end{eqnarray}
We write equation \eqref{wittichequation} as $Q(e^z)=-w''/w-P(e^z)w'/w$. Since $w$ is of finite order, then an application of the lemma on the logarithmic derivative yields $\deg(Q(z))m(r,e^z)\leq \deg(P(z))m(r,e^z)+O(\log r)$, i.e., $[\deg(Q(z))-\deg(P(z))]T(r,e^z)\leq O(\log r)$. Therefore, $\deg(Q(z))\leq \deg(P(z))$ and thus $P(z)$ is nonconstant. Together with equations \eqref{second tumura-cluniesimple4 fube} and \eqref{second tumura-cluniesimple4 fu1be1}, we rewrite equation \eqref{wittichequation} as
\begin{equation}\label{wittichequation re0}
\frac{\sum_{j=0}^k(c^2+2jc+j^2)a_je^{jz}}{\sum_{j=0}^ka_je^{jz}}+\frac{\sum_{j=0}^k(c+j)a_je^{jz}}{\sum_{j=0}^ka_je^{jz}}P(e^z)+Q(e^z)=0.
\end{equation}
Since along a ray $z=re^{i\theta}$ such that $\cos \theta<0$, we have $e^{z}\to 0$ as $r\to\infty$, then by letting $r\to \infty$ along the ray $z=re^{i\theta}$, we obtain from equation \eqref{wittichequation re0} that $c^2+cP(0)+Q(0)=0$. This completes the proof.

\end{proof}

Unlike in sections~\ref{se2: Tumura--Clunie differential equations} and~\ref{se3: an oscillation question of Ishizaki} where Nevanlinna theory plays the central role in proving Theorems~\ref{maintheorem3} and~\ref{maintheorem4}, the proof of Theorem~\ref{maintheorem5} will, however, mainly rely on the \emph{Lommel transformation} for the \emph{generalized Bessel equation}:
\begin{equation}\label{Bank-laine se3simple hobessel1}
x^2y''+xy'+\left(\sum_{-n'}^nd_jx^{j}\right)y=0.
\end{equation}
Recall the \emph{Bessel equation}: $x^2y''+xy'+(x^2-\nu^2)y=0$, where $\nu$ is a nonzero constant. Lommel \cite{Lommel1871} and Pearson \cite{Pearson1880} independently (see also \cite{Watson1944}) studied the following transformation given by:
\begin{equation}\label{Bank-laine se3simple hobessel2}
x=\alpha t^{\beta}, \quad y(x)=t^{\gamma}u(t),
\end{equation}
where $\alpha$, $\beta$ and $\gamma$ are constants and applied to the Bessel equation. By using the above transformation to equation \eqref{Bank-laine se3simple hobessel1} and by computing the derivatives of $x$ and $y$, we get
\begin{equation}\label{Bank-laine se3simple hobessel3}
t^2u''(t)+(2\gamma+1)tu'(t)+\left(\gamma^2+\beta^2\sum_{-n'}^n\alpha^jd_jt^{\beta j}\right)=0.
\end{equation}
A further change of variable such that
\begin{equation}\label{Bank-laine se3simple hobessel4}
t=e^{pz}, \quad f(z)=u(t),
\end{equation}
leads to an equation of the form
\begin{equation}\label{Bank-laine se3simple hobessel5}
f''+2\gamma pf'+p^2\left(\gamma^2+\beta^2\sum_{-n'}^n\alpha^jd_je^{\beta pjz}\right)=0.
\end{equation}
In the case of equation \eqref{Bank-laine se3simple}, by Lommel's transformation we have
\begin{equation}\label{Bank-laine se3simple hobessel5tran1}
x^2y''+xy'-\left(d_1x^{l}+d_2x^{s}+d_3\right)y=0,
\end{equation}
where $d_1$, $d_2$ and $d_3$ are some constants. By comparing the coefficients of equation \eqref{Bank-laine se3simple} and \eqref{Bank-laine se3simple hobessel5}, we deduce that $2\gamma p=0$, $\beta p=1$, $\alpha^{l}d_1=1$, $\alpha^{s}d_2=b_2$ and $d_3=b_3$. Further, for equation \eqref{Bank-laine se3simple hobessel5tran1}, it is well-known that the transformation $y=x^{-1/2}u$ leads to an equation of the form
\begin{equation}\label{Bank-laine se3simple hobessel5tran3}
u''-\left[\frac{1}{\alpha^{l}}x^{l-2}+\frac{b_2}{\alpha^{s}}x^{s-2}+\left(b_3-\frac{1}{4}\right)\frac{1}{x^2}\right]u=0.
\end{equation}
In the case $l=4$, it has been shown by Chiang and Yu~\cite{Chiang2019} that there is a full correspondence between solutions of \eqref{Bank-laine se3simple} such that $\lambda(f)<\infty$ and Liouvillian solutions of \eqref{Bank-laine se3simple hobessel5tran3}. The only possible singular point of equation \eqref{Bank-laine se3simple hobessel5tran3} is $x=0$. Concerning the local solutions around a singular point of a second-order linear differential equation, we have the following elementary Lemma~\ref{lemma  solution0pre0besse3  ahoth0 hhta}; see \cite{Herold1975} or in \cite[Lemma~6.6]{Laine1993}.

\begin{lemma}[\cite{Herold1975,Laine1993}]\label{lemma  solution0pre0besse3  ahoth0 hhta}
Suppose that $h$ is analytic in $|z|<R$, $R>0$, and consider the differential equation
\begin{equation}\label{jahgoahta;aio}
u''+\frac{h(z)}{z^2}u=0
\end{equation}
in the disc $|z|<R$. Let $\rho_1$ and $\rho_2$ be the roots of
\begin{equation}\label{jahgoahta;aiohgia;g}
\rho(\rho-1)+h(0)=0.
\end{equation}
Denote by $D=D(r)$ the slit disc $D:=\{z:|z|<r\}\setminus\{t \ | \ 0\leq t<r\}$. Then
\begin{itemize}

  \item [(1)] if $\rho_1-\rho_2\in \mathbb{Z}\setminus\{0\}$, then equation \eqref{jahgoahta;aio} admits in some slit disc $D=D(r)$, $r\leq R$, two linearly independent solutions $u_1$ and $u_2$ of the form:
\begin{equation}\label{jahgoahta;aiohfa;goie3}
\begin{split}
u_1(z) & =z^{\rho_1}\sum_{i=0}^{\infty}a_iz^i, \quad a_0\not=0, \\
u_2(z) & =u_1(z)d\log z+z^{\rho_2}\sum_{i=0}^{\infty}b_iz^i,
\end{split}
\end{equation}
where $d=0$ or $d=1$;
  \item [(2)] if $\rho_1-\rho_2\not\in \mathbb{Z}$, then equation \eqref{jahgoahta;aio} admits in some slit disc $D=D(r)$, $r\leq R$, two linearly independent solutions $u_1$ and $u_2$ of the form:
\begin{equation}\label{jahgoahta;aiohfa;goie2}
\begin{split}
u_1(z) & =z^{\rho_1}\sum_{i=0}^{\infty}a_iz^i, \quad a_0\not=0, \\
u_2(z) & =z^{\rho_2}\sum_{i=0}^{\infty}b_iz^i, \quad b_0\not=0;
\end{split}
\end{equation}
  \item [(3)] if $\rho_1-\rho_2=0$, then equation \eqref{jahgoahta;aio} admits in some slit disc $D=D(r)$, $r\leq R$, two linearly independent solutions $u_1$ and $u_2$ of the form:
\begin{equation}\label{jahgoahta;aiohfaaj;goie1}
\begin{split}
u_1(z) & =z^{\rho_1}\sum_{i=0}^{\infty}a_iz^i, \quad a_0\not=0, \\
u_2(z) & =u_1(z)\log z+z^{\rho_2}\sum_{i=0}^{\infty}b_iz^i.
\end{split}
\end{equation}
\end{itemize}

\end{lemma}

For the solution $u_2$ in~\eqref{jahgoahta;aiohfa;goie3}, if $d=0$, then from the proof of \cite[Lemma~6.6]{Laine1993} we know that $b_0\not=0$.

Now, by elementary theory of ordinary differential equation (see \cite{Herold1975}), Lemma~\ref{lemma  solution0pre0besse3  ahoth0 hhta} shows that equation~\eqref{Bank-laine se3simple hobessel5tran3} admits two linearly independent solutions $u_1$ and $u_2$ in the broken plane $\mathbb{C}^{-}=\mathbb{C}\setminus\{x \ | \ 0\leq x<\infty\}$. When $p=1$ in \eqref{Bank-laine se3simple hobessel4}, by the Lommel transformation and analytic continuation principle, the general solution of \eqref{Bank-laine se3simple} is thus given by
\begin{equation}\label{Bank-laine se3simple hobessel5a;goija;oog0}
f=(\alpha e^{z})^{-1/2}[E_1u_1(\alpha e^{z})+E_2u_2(\alpha e^{z})],
\end{equation}
where $E_1$ and $E_2$ are two arbitrary constants. Note that the above solution is independent from the choice of the branches of $u_1$ and $u_2$ in Lemma~\ref{lemma  solution0pre0besse3  ahoth0 hhta}. This is the key observation for the proof of Theorem~\ref{maintheorem5}.

Now we begin to prove Theorem~\ref{maintheorem5}.

\renewcommand{\proofname}{Proof of Theorem~\ref{maintheorem5}.}
\begin{proof}

We first suppose that $f$ is a nontrivial solution such that $\lambda(f)<\infty$ of equation~\eqref{Bank-laine se3simple} and use the expressions in \eqref{bank-laine01hf1} and \eqref{bank-laine01hf2} to write $f(z)=\Psi(x)=x^{c}\psi(x)e^{\chi(x)}$, where $x=e^{z/h}$, $h=1$ or $h=2$. From the proof of \cite[Theorem~1.2]{Chiang2019}, we know that in the broken plane $\mathbb{C}^{-}$ equation \eqref{Bank-laine se3simple hobessel5tran3} admits a solution of the form
\begin{equation}\label{Bank-laine sjhgfore}
u=\exp\left(\int \omega dx\right),
\end{equation}
where $\omega:=\chi'+\psi'/\psi+(2hc+1)/(2x)$ is rational function in the complex plane~$\mathbb{C}$.
By using Kovacic's algorithm in \cite{Kovacic1986} and giving the same discussions as in the proof of \cite[Theorem~3.1]{Chiang2019} to equation \eqref{Bank-laine se3simple hobessel5tran3} for the two cases $b_3\not=1/4$ and $b_3=1/4$, respectively, we conclude that $l$ must be even. It follows that $f$ has the form in \eqref{bank-laine01hf2} and thus $p=1$ in \eqref{Bank-laine se3simple hobessel4} and $\alpha=1$ in~\eqref{Bank-laine se3simple hobessel2}. We write $f=\kappa e^{h}$, where $\kappa$ and $h'$ are both polynomials in $\zeta=e^z$ such that $\kappa(0)\not=0$. We may also write $f=\kappa_ce^{h_c}$, where $\kappa_c=\kappa e^{cz}$ and $h_c'=h'-c$. Then, denoting $g_c=h_c'$, we have from \eqref{Bank-laine se3simple} that
\begin{equation}\label{second tumura-clunieja;gij0jjjh1}
g_c^2+g'_c+2\frac{\kappa'_c}{\kappa_c}g_c+\frac{\kappa''_c}{\kappa_c}=e^{lz}+b_2e^{sz}+b_3.
\end{equation}
By Lemma~\ref{lemma  solution0pre} and the expression in \eqref{bank-laine01hf2}, we see that the constant $c$ in \eqref{bank-laine01hf2} satisfies $c^2=b_3$.

Now, for the solution in \eqref{Bank-laine se3simple hobessel5a;goija;oog0}, by Lemma~\ref{lemma  solution0pre0besse3  ahoth0 hhta} we have~$\rho_1+\rho_2=1$ and $\rho_1\rho_2=1/4-b_3$, which yield $(\rho_1-\rho_2)^2=4b_3$. Then $\rho_1-\rho_2=-2c$ and it follows that $\rho_1=(1-2c)/2$ and $\rho_2=(1+2c)/2$. Thus the solutions in \eqref{Bank-laine se3simple hobessel5a;goija;oog0} can be written as
\begin{equation}\label{Bank-laine se3simple hobessel5a;goija;oog0jh}
\begin{split}
f=(e^{z})^{c}\left[\left(E_1+E_2d\log e^z\right)(e^{z})^{-2c}\sum_{j=0}^{\infty}a_j(e^{z})^{j}+E_2\sum_{j=0}^{\infty}b_j(e^{z})^{j}\right],
\end{split}
\end{equation}
when $\rho_1-\rho_2\not=0$, or
\begin{equation}\label{Bank-laine se3simple hobessel5a;goija;oog0jh1}
\begin{split}
f=\left(E_1+E_2\log e^z\right)\sum_{j=0}^{\infty}a_j(e^{z})^{j}+E_2\sum_{j=0}^{\infty}b_j(e^{z})^{j},
\end{split}
\end{equation}
when $\rho_1-\rho_2=0$. Note that $d=0$ in \eqref{Bank-laine se3simple hobessel5a;goija;oog0jh} when $\rho_1-\rho_2$ is not an integer. On the other hand, for the solution $f=\kappa e^{h}$, we may write the expression in \eqref{bank-laine01hf2} in the form $f=e^{cz}\sum_{j=0}^{\infty}d_je^{jz}$. By comparing this series with the one in \eqref{Bank-laine se3simple hobessel5a;goija;oog0jh} or in \eqref{Bank-laine se3simple hobessel5a;goija;oog0jh1}, we conclude that the logarithmic term in $u_2$ does not occur. This implies that $d=0$ or $E_2=0$ in \eqref{Bank-laine se3simple hobessel5a;goija;oog0jh} and $E_2=0$ in \eqref{Bank-laine se3simple hobessel5a;goija;oog0jh1} when $f=\kappa e^{h}$.

With these preparations, we now suppose that $f_1$ and $f_2$ are two linearly independent solutions of~\eqref{Bank-laine se3simple} such that $\max\{\lambda(f_1),\lambda(f_2)\}<\infty$. We write $l=2(m+1)$ for some integer $m\geq 0$ and also write $s=2(m+1)-t$ for some integer $t\geq 1$.

Let $q$ be the smallest integer such that $s/l\leq [2(q+1)-1]/[2(q+1)]$. Since $l$ and $s$ are relatively prime, we see that the equality holds only when $q=m$. For each of $f_1$ and $f_2$, denoted by $f$, we may write $f=\kappa_ce^{h_c}$, where $\kappa_c=\kappa e^{cz}$. Then, denoting $g_c=h_c'$, we have equation \eqref{second tumura-clunieja;gij0jjjh1}. By Theorem~\ref{maintheorem3}, $g_c=h_c'=\sum_{j=0}^{q}c_je^{(m+1-jt)z}$, where $c_0,c_1,\cdots,c_q$ are constants such that $c_0^2=1$ and $c_1,\cdots,c_q$ satisfy certain relations. In both of the two cases where $q=1$ and $q\geq 2$, we have $2c_0c_1=b_2$ and, by simple computations, that,
\begin{equation}\label{coeff1}
\begin{split}
c_j=\frac{t_jc_1^{j}}{(-2c_0)^{j-1}}, \quad  j=1,\cdots,q,
\end{split}
\end{equation}
where $t_j$ are positive integers such that $t_1<\cdots<t_q$, and further that
\begin{equation}\label{coeff2}
\begin{split}
C_{q+j}=\sum_{\substack{j_0+\cdots+j_q=2,\\j_1+\cdots+qj_q=q+j}}\frac{2}{j_0!\cdots j_q!}c_0^{j_0}\cdots c_q^{j_q}=\frac{T_{q+j}c_1^{q+j}}{(-2c_0)^{q+j-2}}, \quad  j=1,\cdots,q,
\end{split}
\end{equation}
where $T_{q+j}$ are positive integers such that $T_{q+1}<\cdots<T_{2q}$. By substituting $g_c=\sum_{j=0}^{q}c_je^{(m+1-jt)z}$ into \eqref{second tumura-clunieja;gij0jjjh1} together with Theorem~\ref{maintheorem3}, we get
\begin{equation}\label{second tumura-cluniesimplej;aojta;io2jyrqbe}
\begin{split}
\frac{\kappa_c''}{\kappa_c}+2c_0e^{(m+1)z}\frac{\kappa_c'}{\kappa_c}+c_0(m+1)e^{(m+1)z}-b_2e^{sz}-b_3=0,
\end{split}
\end{equation}
when $q=0$, and
\begin{equation}\label{second tumura-cluniesimplej;aojta;io2jyrq}
\begin{split}
&\frac{\kappa_c''}{\kappa_c}+2\left(\sum_{j=0}^{q}c_je^{(m+1-jt)z}\right)\frac{\kappa_c'}{\kappa_c}-b_3\\
&+\sum_{j=0}^{q}\left[C_{k_j}e^{[m+1-(q+1)t]z}+(m+1-jt)c_j\right]e^{(m+1-jt)z}=0,
\end{split}
\end{equation}
when $q\geq1$, where $C_{k_j}=C_{q+1+j}$ and $C_{2q+1}=0$. By substituting equations \eqref{second tumura-cluniesimple4 fube} and \eqref{second tumura-cluniesimple4 fu1be1} for $\kappa_{c}=e^{cz}(\sum_{i=1}^ka_ie^{iz})$, $a_0a_k\not=0$, into \eqref{second tumura-cluniesimplej;aojta;io2jyrqbe} or \eqref{second tumura-cluniesimplej;aojta;io2jyrq} and noting that $b_3=c^2$, we finally get
\begin{equation}\label{second tumura-clunieja;gij1a11a}
\begin{split}
\sum_{i=0}^k\left[(2ic+i^2)a_ie^{iz}+c_0(2c+2i+m+1)a_ie^{(m+1+i)z}-b_2a_ie^{(i+s)z}\right]=0,
\end{split}
\end{equation}
when $q=0$, and
\begin{equation}\label{second tumura-cluniesimple4 fu2jyq}
\begin{split}
&\sum_{i=0}^k(2ic+i^2)a_ie^{iz}+2\left(\sum_{i=0}^k(c+i)a_ie^{iz}\right)\left(\sum_{j=0}^{q}c_je^{(m+1-jt)z}\right)\\
&+\left(\sum_{i=0}^ka_ie^{iz}\right)\left\{\sum_{j=0}^{q}[C_{k_j}e^{[m+1-(q+1)t]z}+(m+1-jt)c_j]e^{(m+1-jt)z}\right\}=0,
\end{split}
\end{equation}
when $q\geq 1$. Note that the inequality $(2q-1)/(2q)<s/l\leq [2(q+1)-1]/[2(q+1)]$ implies $qt<m+1\leq (q+1)t$, where the equality holds when $q=m$. The left-hand side of equations \eqref{second tumura-clunieja;gij1a11a} and \eqref{second tumura-cluniesimple4 fu2jyq} are polynomials in $e^{z}$ of degree $k+m+1$ and thus all coefficients of these two polynomials vanish. When $s=2(m+1)-t<2m+1$, we have $q<m$ and $t\geq 2$. By looking at the highest-degree term in the resulting polynomial and noting that $a_k\not=0$, we find
\begin{equation}\label{second tumura-cluniesimple4 fu2jyq2}
\begin{split}
m+2c+2k+1=0.
\end{split}
\end{equation}
Similarly, when $s=2m+1$, we have $q=m$ and $t=1$ and find
\begin{equation}\label{second tumura-cluniesimple4 fu2jyq1}
\begin{split}
C_{m+1}+(m+2c+2k+1)c_0=0.
\end{split}
\end{equation}
Let $c_{+}$ or $c_{-}$ be any square root of $b_3$. We may write $f_1=\kappa_{1c}e^{h_{1c}}$ and $f_2=\kappa_{2c}e^{h_{2c}}$, where $\kappa_{1c}=\kappa_1e^{c_{+}z}$ and $\kappa_{2c}=\kappa_2e^{c_{-}z}$ and $\kappa_{1}$ and $\kappa_{2}$ are two polynomials in $e^z$ of degrees $k_1$ and $k_2$, respectively. Moreover, $h_{1c}$ satisfies $h_{1c}'=\sum_{j=0}^{q}c_je^{(m+1-jt)z}$. Since $c_0^2=1$ and since $2c_0c_1=b_2$ when $q\ge1$, we easily deduce from \eqref{coeff1} that $h_{2c}=\pm h_{1c}$. Recall the elementary \emph{Wronskian determinant}: $f_1'f_2-f_1f_2'=D$, where~$D$ is a nonzero constant (see \cite{Laine1993}). Then we have $(f_1/f_2)'=D/f_2^2$. If $h_{2c}=h_{1c}$, then $f_1/f_2$ is of finite order while $f_1^2$ is of infinite order, a contradiction. Therefore, $h_{2c}=-h_{1c}$. We may suppose that $c_{+}=c$. When $s<2m+1$, from \eqref{second tumura-cluniesimple4 fu2jyq2} we deduce that $c_{+}=c_{-}=c$, which implies that $-2c=m+1+2k_1$ is a positive integer and hence $k_1=k_2$; when $s=2m+1$, using \eqref{coeff2} and the relation $2c_0c_1=b_2$ we deduce from \eqref{second tumura-cluniesimple4 fu2jyq1} that $c_{+}+c_{-}+m+1+k_1+k_2=0$, which implies that $c_{+}=c_{-}=c$ and hence $-2c=m+1+k_1+k_2$ is a positive integer. Now $\rho_1-\rho_2=-2c$ is a positive integer. Together with \eqref{Bank-laine se3simple hobessel5a;goija;oog0} and previous preparations, we conclude that equation \eqref{Bank-laine se3simple hobessel5tran3} admits in the broken plane $\mathbb{C}^{-}$ two linearly independent solutions of the form $u_1=x^{\rho_1}v_1$ and $u_2=x^{\rho_2}v_2$, where $v_1$ and $v_2$ are two entire functions such that $v_1(0)\not=0$ and $v_2(0)\not=0$, so that
\begin{equation}\label{coeff1uy1}
\begin{split}
f_1&=x^{c}\kappa_{11}e^{h_{11}}=x^{-1/2}\left(D_1x^{\rho_1}v_1+D_2x^{\rho_2}v_2\right),\\
f_2&=x^{c}\kappa_{12}e^{-h_{11}}=x^{-1/2}\left(D_3x^{\rho_1}v_1+D_4x^{\rho_2}v_2\right),
\end{split}
\end{equation}
where $D_j$ are constants, $h_{11}=\sum_{j=0}^{q}\frac{c_j}{m+1-jt}x^{m+1-jt}$, and $\kappa_{11}$ and $\kappa_{12}$ are two polynomials of degrees $k_1$ and $k_2$, respectively, such that $\kappa_{11}(0)\not=0$ and $\kappa_{12}(0)\not=0$. Noting $\rho_1=(1-2c)/2$ and $\rho_2=(1+2c)/2$, we see from \eqref{Bank-laine se3simple hobessel5a;goija;oog0jh} that $D_2D_4\not=0$. Obviously, $E:=D_1D_4-D_2D_3\not=0$. From equation \eqref{coeff1uy1} we get
\begin{equation}\label{coeff1uy2}
\begin{split}
u_1=x^{\rho_1}v_1=\frac{1}{E}x^{1/2}x^{c}\left(D_4\kappa_{11}e^{h_{11}}-D_2\kappa_{12}e^{-h_{11}}\right).
\end{split}
\end{equation}
Since $v_1$ is an entire function with $v_1(0)\not=0$, we see from \eqref{coeff1uy2} that the function $w:=D_4\kappa_{11}e^{2h_{11}}-D_2\kappa_{12}$ has a zero of order $-2c$ at the point $z=0$ and so $w(0)=w'(0)=\cdots=w^{(-2c-1)}(0)=0$. Denote
\begin{equation}\label{coeff1uy3}
\begin{split}
\kappa_{11}&=a_{1,0}+a_{1,1}x+\cdots+a_{1,k_1}x^{k_1},\\
\kappa_{12}&=a_{2,0}+a_{2,1}x+\cdots+a_{2,k_2}x^{k_2},
\end{split}
\end{equation}
where $a_{1,0},a_{2,0},a_{1,k_1},a_{1,k_2}\not=0$. $w(0)=0$ implies that $D_4a_{1,0}=D_2a_{2,0}$. Supposing that $a_{1,0}=a_{2,0}=1$, we have $D_4=D_2$. Below we consider the case when $m\geq 1$.

Consider first the case when $s/l<1/2$. Since $m\geq 1$, by Theorem~\ref{maintheorem1} and \eqref{second tumura-cluniesimple4 fu2jyq2}, we see that $k_1=k_2\geq 1$. Now $w^{(m+1)}(0)=0$ implies that $a_{1,m+1}+2(m!)c_0=a_{2,m+1}$. Here $a_{1,m+1}=0$ if $m+1>k_1$ an so is for $a_{2,m+1}$. Obviously, $m+1\leq k_1$. For each of $f_1$ and $f_2$, denoted by $f$, we may write $f=\kappa_ce^{h_c}$. Then we have equation~\eqref{second tumura-clunieja;gij1a11a}. Note that $1\leq s=2(m+1)-t\leq m$. The left-hand side of equation~\eqref{second tumura-clunieja;gij1a11a} is a polynomial in $e^{z}$ of degree $m+1+k$ and thus all coefficients of this polynomial vanish. Denoting $a_{-m-1}=\cdots=a_{-2}=a_{-1}=0$ and $a_{k+1}=\cdots=a_{k+m}=0$, we obtain from equation \eqref{second tumura-clunieja;gij1a11a} that $c_0(2c+2k+m+1)a_k=0$, which yields \eqref{second tumura-cluniesimple4 fu2jyq2}, and
\begin{equation}\label{second tumura-clunieja;gij1a14sagae}
\begin{split}
c_0(2c+2i-m-1)a_{i-m-1}&=b_2a_{i-s}-(2ic+i^2)a_{i},  \quad i=0,\cdots,k+m.
\end{split}
\end{equation}
By substituting $2c=-2k-m-1$ into the equations in \eqref{second tumura-clunieja;gij1a14sagae} we get
\begin{equation}\label{second tumura-clunieja;gij1a15afgg}
\begin{split}
2c_0(k-i+m+1)a_{i-m-1}=-b_2a_{i-s}+(2ic+i^2)a_i,  \quad i=0,\cdots,k+m.
\end{split}
\end{equation}
Since $1\leq s\leq m$, by letting $i=1,2,\cdots,m$, we see that $a_1/a_0=K_1$, $\cdots$, $a_m/a_0=K_m$ for some constants $K_1,\cdots,K_m$ independent from $c_0$. Then by letting $i=m+1$ in \eqref{second tumura-clunieja;gij1a15afgg} together with the relation $a_{1,m+1}+2(m!)c_0=a_{2,m+1}$, we get $2k_1/(m+1)(2c+m+1)+2(m!)=-2k_2/(m+1)(2c+m+1)$. Since $k_1=k_2$ and $-2c=2k_1+m+1$, we get $(m+1)!=1$, which is impossible when $m\geq 1$.

Consider next the case when $s/l>1/2$ and $s<2m+1$. Now $qt<m+1<(q+1)t$ for some integer $1\leq q<m$. Recalling that $s=2(m+1)-t$ and $l$ and $s$ are relatively prime, we see that $t\geq3$ is an odd integer. Denoting each of $f_1$ and $f_2$ by $f$, we may write $f=\kappa_ce^{h_c}$. Then we have equation \eqref{second tumura-cluniesimple4 fu2jyq}. Denote $M=m+1-qt$ for simplicity. Since $-2c=2k+m+1$ and $C_{2q+1}=0$ and $C_{2q}=c_q^2$, then by
looking at the coefficient of the term $e^{Mz}$ on the left-hand side of equation \eqref{second tumura-cluniesimple4 fu2jyq}, we find $M(2c+M)a_M+2ca_0c_q+Ma_0c_q=0$, which gives $a_0c_q+Ma_M=0$.
Recall that the function $w:=D_4\kappa_{11}e^{2h_{11}}-D_2\kappa_{12}$ has a zero of order $-2c$ at the point $z=0$. Then $w^{(M)}(0)=0$ implies that $a_{1,M}+2(M-1)!c_q=a_{2,M}$. Using equation \eqref{coeff1} together with $2c_0c_1=b_2$ and $a_0c_q+Ma_M=0$, we get $M!=1$, which implies that $M=1$. It follows that $m=qt$. By looking at the coefficient of the term $e^{2z}$ on the left-hand side of equation \eqref{second tumura-cluniesimple4 fu2jyq}, we find $2(2c+2)a_2+2(c+1)a_1c_q+a_0c_q^2+a_1c_q=0$, which together with $a_1=-a_0c_q$ yields $2a_2=a_0c_q^2$. Then by looking at the coefficient of the term $e^{3z}$ on the left-hand side of equation \eqref{second tumura-cluniesimple4 fu2jyq}, we find $3(2c+3)a_3+2(c+2)a_2c_q+a_2c_q+a_1c_q^2=0$, i.e., $(2c+3)(6a_3+c_q^3)=0$ and thus $6a_3+c_q^3=0$. Now $w^{(3)}(0)=0$ implies that $a_{1,3}+6a_{1,2}c_q+12a_{1,1}c_q^2+8c_q^3=a_{2,3}$, which together with the relation $a_1=-a_0c_q$ and $2a_2=a_0c_q^2$ gives $a_{1,3}-c_q^3=a_{2,3}$. Then using equation \eqref{coeff1} together with $2c_0c_1=b_2$ and $c_q^3+6a_3=0$, we get $c_q^3=0$, a contradiction to \eqref{coeff1}.

Finally, we consider the case when $s=2m+1$. Recall that $q=m$ and $t=1$ in \eqref{second tumura-cluniesimple4 fu2jyq}. In this case, if $k_1=k_2$, then using \eqref{coeff2} and the relation $2c_0c_1=b_2$ we get from \eqref{second tumura-cluniesimple4 fu2jyq1} that $C_{m+1}=0$, a contradiction. Therefore, without loss of generality, we may suppose that $k_1>k_2\geq 0$. If $k_2=0$, then by Theorem~\ref{maintheorem1} have $2c+1=0$, which is impossible since $-2c=m+1+k_1+k_2$. Therefore, $k_2>0$. Note that $C_{2m}=c_m^2$. Since $-2c=m+1+k_1+k_2$, then by looking at the coefficient of the terms $e^z$ and $e^{2z}$ in equation \eqref{second tumura-cluniesimple4 fu2jyq}, respectively, we find $a_{1}+a_{0}c_m=0$ and $2(2c+2)a_2+(2c+3)c_ma_1+[c_m^2+(2c+2)c_{m-1}]a_0=0$ and so $2a_2-c_m^2a_0+c_{m-1}a_0=0$. Recall that the function $w:=D_4\kappa_{11}e^{2h_{11}}-D_2\kappa_{12}$ has a zero of order $-2c$ at the point $z=0$. Now $w''(0)=0$ implies that $a_{1,2}+4a_{1,1}c_m+4c_{m-1}+4c_m^2=a_{2,2}$. Using equation \eqref{coeff1} together with $2c_0c_1=b_2$, we get $-c_{m-1}/2+4c_{m-1}=c_{m-1}/2$, which yields $c_{m-1}=0$, a contradiction to \eqref{coeff1}.

From the above reasoning, we conclude that $l=2$. We complete the proof.

\end{proof}

In the rest of this section, we use Theorem~\ref{maintheorem3} to determine precisely all nontrivial solutions such that $\lambda(f)<\infty$ of equation \eqref{Bank-laine se3simple} for the case $l=2$ and $l=4$.

\begin{theorem}\label{maintheorem6}
Let $b_1$, $b_2$ and $b_3$ be constants such that $b_1b_2\not=0$ and $s$ and $l$ be relatively prime integers such that $1\leq s<l\leq 4$. Suppose that~\eqref{Bank-laine se3simple} admits a nontrivial solution~$f$ such that $\lambda(f)<\infty$. Then
\begin{itemize}

\item [(1)]
if $s=1$ and $l=2$, then $f=\kappa e^{h}$, $\kappa=\sum_{i=-1}^ka_ie^{iz}$ and $h=c_0e^{z}+cz$, where $k\geq 0$ is an integer, $c_0$ and $c$ are constants such that $c_0^2=1$, $2c_0(c+k)+c_0=b_2$ and $c^2=b_3$, and $a_{-1}$, $a_0$, $\cdots$, $a_k$ are constants such that $a_0a_k\not=0$, $a_{-1}=0$ and
\begin{equation}\label{recuree1}
\begin{split}
2c_0(k+1-i)a_{i-1}=(2ic+i^2)a_i, \quad i=0,1,\cdots,k;
\end{split}
\end{equation}

\item [(2)]
if $s=1$ and $l=4$, then $f=\kappa e^{h}$, $\kappa=\sum_{i=-2}^{k+1}a_ie^{iz}$ and $h=(c_0/2)e^{2z}+cz$, where $k\geq 1$ is an integer, $c_0$ and $c$ are constants such that $c_0^2=1$, $2c+2k+2=0$ and $c^2=b_3$, and $a_{-2}$, $a_{-1}$, $a_0$, $\cdots$, $a_{k+1}$ are constants such that $a_0a_k\not=0$, $a_{-2}=a_{-1}=a_{k+1}=0$ and
\begin{equation}\label{recuree2}
\begin{split}
2c_0(k-i+2)a_{i-2}=-b_2a_{i-1}+(2ic+i^2)a_i,  \quad i=0,1,\cdots,k+1;
\end{split}
\end{equation}

\item [(3)]
if $s=3$ and $l=4$, then $f=\kappa e^{h}$, $\kappa=\sum_{i=-2}^{k+1}a_ie^{iz}$ and $h=(c_0/2)e^{2z}+c_1e^{z}+cz$, where $k\geq 0$ is an integer, $c_0$, $c_1$ and $c$ are constants such that $c_0^2=1$, $2c_0c_1=b_2$, $c^2=b_3$ and $c_1^2+(2+2c+2k)c_0=0$, and $a_{-2}$, $a_{-1}$, $a_0$, $\cdots$, $a_{k+1}$ are constants such that $a_0a_k\not=0$, $a_{-2}=a_{-1}=a_{k+1}=0$ and
\begin{equation}\label{recuree3}
\begin{split}
(2k-2i+4)c_0a_{i-2}=(2c+2i-1)c_1a_{i-1}+(2ic+i^2)a_{i}, \quad  i=0,\cdots,k+1.
\end{split}
\end{equation}

\end{itemize}

\end{theorem}

For the convenience to write the recursive formulas in \eqref{recuree1}-\eqref{recuree3}, we have introduced some extra coefficients $a_{-2}$, $a_{-1}$, $a_{k+1}$, which are all equal to~$0$.

When $s=1$ and $l=4$, since $a_0\not=0$ and $a_{k}\not=0$ and $2c+2k+2=0$, the recursive formulas in \eqref{recuree2} yield a polynomial equation $P(b_2)=0$ with respect to $b_2$ with coefficients formulated in terms of $c_0$ and $k$. For example, when $k=1$, we have $0=-b_2a_{0}+(2c+1)a_1$ and $2c_0a_{0}=-b_2a_{1}$, which together with the equation $2c+4=0$ yield $b_2^2-6c_0=0$, etc.

When $s=3$ and $l=4$, if $2c+1\not=0$, then we may solve from the first $k+1$ equations in \eqref{recuree3} that $a_{k-1}=P(c)a_k$ for some polynomial $P(c)$ with respect to $c$ with coefficients formulated in terms of $c_0$ and $c_1$. By combining this equation with the equation $2c_0a_{k-1}=(2c+k)c_1a_{k}$ together with the relation $c_1^2+(2+2c+2k)c_0=0$ we may obtain a polynomial equation $P(t)=0$ with respect to $t=2c$ with coefficients independent from $c_0$ and $c_1$. For example, when $k=1$, we have $0=(2c+1)(c_1a_{0}+a_{1})$ and $2c_0a_{0}=(2c+3)c_1a_{1}$, which yield $P(t)=(t+2)(t+5)=0$ and thus $2c=-2$ or $2c=-5$, etc.


\renewcommand{\proofname}{Proof of Theorem~\ref{maintheorem6}.}
\begin{proof}
Suppose that $f$ is a nontrivial solution such that $\lambda(f)<\infty$ of equation~\eqref{Bank-laine se3simple}. 
Following the proof of Theorem~\ref{maintheorem5}, we may write $f=\kappa_c e^{h_c}$, where $\kappa_c=\kappa e^{cz}$ and $g_{c}=h_{c}'$ and then from \eqref{Bank-laine se3simple} we get equation \eqref{second tumura-clunieja;gij0jjjh1}. Below we consider three cases: (1) $s=1$ and $l=2$; (2) $s=1$ and $l=4$; (3) $s=3$ and $l=4$.

For the first two cases $s=1$ and $l=2$ or $s=1$ and $s=4$, we have $\kappa_c=e^{cz}(\sum_{i=0}^ka_ie^{iz})$ and $h_c=[c_0/(m+1)]e^{(m+1)z}$, where $m=0$ or $m=1$ and $c_0$ is a constant such that $c_0^2=1$. Moreover, if $m=1$, then by Theorem~\ref{maintheorem1} we see $k\geq1$. From the proof of Theorem~\ref{maintheorem5}, we have equations \eqref{second tumura-cluniesimplej;aojta;io2jyrqbe} and \eqref{second tumura-clunieja;gij1a11a} with $s=1$.
When $l=2$, since the left-hand side of equation~\eqref{second tumura-clunieja;gij1a11a} is a polynomial in $e^{z}$ of degree $1+k$, all coefficients of this polynomial vanish. Therefore, denoting $a_{-1}=0$, we obtain from equation \eqref{second tumura-clunieja;gij1a11a} that $[2c_0(c+k)+c_0-b_2]a_k=0$ and
\begin{equation}\label{second tumura-clunieja;gij1a12}
\begin{split}
[2c_0(c+i-1)+c_0-b_2]a_{i-1}+(2ic+i^2)a_i=0,  \quad i=0,1,\cdots,k.
\end{split}
\end{equation}
Since $a_k\not=0$, we have $2c_0(c+k)+c_0-b_2=0$ and then obtain the recursive formulas in \eqref{recuree1} by substituting $2c_0c+c_0-b_2=-2c_0k$ into the equations in \eqref{second tumura-clunieja;gij1a12}. When $m=1$, we have the recursive formulas in \eqref{second tumura-clunieja;gij1a15afgg} with $s=1$. Denoting $a_{-2}=a_{-1}=0$ and $a_{k+1}=0$, we have the recursive formulas in \eqref{recuree2}.

When $s=3$ and $l=4$, we have $\kappa_c=e^{cz}(\sum_{i=0}^ka_ie^{iz})$ and $h_c=(c_0/2)e^{2z}+c_1e^{z}$, where $c_0$, $c_1$ are two constants such that $c_0^{2}=1$, $2c_0c_1=b_2$. From the proof of Theorem~\ref{maintheorem5}, we get equation \eqref{second tumura-cluniesimple4 fu2jyq} with $q=m=1$. Similarly as in previous cases, denoting $a_{-2}=a_{-1}=a_{k+1}=0$, we finally get the recursive formulas in \eqref{recuree3}. We omit those details.

\end{proof}

By Theorem~\ref{maintheorem6}, we may give a different formulation from the results in~\cite[Theorem~1.6]{Chiang2006}.

\begin{corollary}\label{corollary1}
Let $s=1$ and $l=2$. Then equation \eqref{Bank-laine se3simple} admits two linearly independent solutions $f_1$ and $f_2$ such that $\max\{\lambda(f_1),\lambda(f_2)\}<\infty$ if and only if there are two distinct nonnegative integers $k_1,k_2$ such that $b_2=\pm(k_1-k_2)$ and $4b_3=(k_1+k_2+1)^2$. In particular, it is possible that $\min\{\lambda(f_1),\lambda(f_2)\}=0$.
\end{corollary}

\renewcommand{\proofname}{Proof of Corollary~\ref{corollary1}.}
\begin{proof}

Let $f_1$ and $f_2$ be two linearly independent solutions of equation \eqref{Bank-laine se3simple} such that $\max\{\lambda(f_1),\lambda(f_2)\}<\infty$. Let $c_{+}$ or $c_{-}$ be any square-root of $b_3$. By Theorem~\ref{maintheorem6}, we may write $f_1=\kappa_{1}e^{h_{1}}$ and $f_2=\kappa_{2}e^{h_{2}}$, where $h_{1}=c_0e^{z}+c_{+}z$ and $c_0$ is a constant such that $c_0^2=1$, $h_{2}=\pm c_0e^{z}+c_{-}z$, $\kappa_{1}$ and $\kappa_2$ are two polynomials in $e^z$ of degrees $k_1$ and $k_2$, respectively. From the proof of Theorem~\ref{maintheorem5} we know that $h_{2}=-c_0e^{z}+c_{-}z$. Since $2c_0(c_{+}+k_1)+c_0=-2c_0(c_{-}+k_2)-c_0=b_2$, we see that $c_{+}=c_{-}$ for otherwise we have $1+k_1+k_2=0$, which is impossible. Letting $c_{+}=c_{-}=c$, then we have $2c+k_{1}+k_{2}+1=0$ and it follows that $b_2=c_0(k_1-k_2)$. Since $b_2\not=0$ and $b_3=c^2$, we have $k_1\not=k_2$ and $4b_3=(k_1+k_2+1)^2$.

Conversely, we let $k_1$ and $k_2$ be two nonnegative integers such that $2c+k_{1}+k_{2}+1=0$, where $c$ satisfies $c^2=b_3$. We first define $f_1=\kappa_1e^{h_1}$, where $\kappa_1=\sum_{i=-1}^{k_1}a_ie^{iz}$, $h_{1}=c_0e^{z}+cz$, $k_1\geq 0$ is an integer, $c_0$ satisfies $c_0^2=1$ and $c_0[2(c+k_1)+1]=b_2$, and $a_{-1}$, $a_0$, $\cdots$, $a_k$ are constants such that $a_{-1}=0$ and
\begin{equation}\label{recuree1pm}
\begin{split}
2c_0(k_1+1-i)a_{i-1}=(2ic+i^2)a_i, \quad i=0,1,\cdots,k_1.
\end{split}
\end{equation}
Also, we define $f_2=\kappa_2 e^{h_2}$, where $\kappa_2=\sum_{i=-1}^{k_2}\hat{a}_ie^{iz}$ and $h_2=-c_0e^{z}+cz$, $k_2\geq 0$ is an integer, $c_0$ satisfies $c_0^2=1$ and $-c_0[2(c+k_2)+1]=b_2$, and $\hat{a}_{-1}$, $\hat{a}_0$, $\cdots$, $\hat{a}_k$ are constants such that $\hat{a}_{-1}=0$ and
\begin{equation}\label{recuree1pm1}
\begin{split}
-2c_0(k_2+1-i)\hat{a}_{i-1}=(2ic+i^2)\hat{a}_i, \quad i=0,1,\cdots,k_2.
\end{split}
\end{equation}
Then by Theorem~\ref{maintheorem6} we see that $f_1$ and $f_2$ are two linearly independent solutions of \eqref{Bank-laine se3simple} such that $\max\{\lambda(f_1),\lambda(f_2)\}<\infty$. Obviously, we may choose one of $k_1$ and $k_2$ to be zero and thus $\min\{\lambda(f_1),\lambda(f_2)\}=0$. We complete the proof.


\end{proof}

\section{Concluding remarks}\label{se5: concluding remarks}

The oscillation of certain second-order linear differential equation \eqref{bank-laine0} are investigated in this paper. If equation \eqref{bank-laine0} with $A(z)$ being a linear combination of two exponential type functions admits a nontrivial solution such that $\lambda(f)<\infty$, by Hadamard's factorization theorem we obtain a Tumura--Clunie type differential equation with coefficients being combinations of functions in $\mathcal{S}$. In section~\ref{se2: Tumura--Clunie differential equations}, we give the form of entire solutions of the Tumura--Clunie type differential equations. As an application, in
section~\ref{se3: an oscillation question of Ishizaki} we give a partial answer to an oscillation question concerning equation \eqref{Bank-laine se3} proposed by Ishizaki~\cite{Ishizaki19970}.
In section~\ref{se4: equation with periodic coefficients}, we consider equation \eqref{bank-laine0} for the case $A(z)=e^{lz}+b_2e^{sz}+b_3$, where $l$ and $s$ are two relatively prime integers and $b_2,b_3$ are constants such that $b_2\not=0$. The general form of solutions such that $\lambda(f)<\infty$ are known. If there are two linearly independent such solutions, we prove that the only possible case is when $l=2$.


By doing straightforward computations, we precisely characterize all solutions such that $\lambda(f)<\infty$ of equation \eqref{Bank-laine se3simple} for the two cases $l=2$ and $l=4$. Unfortunately, we are unable to include or exclude other possibilities. Although, by using Theorems~\ref{maintheorem4} and~\ref{maintheorem5} together with Lemma~\ref{lemma  solution0pre}, when $l\not=2,4$, we may also obtain some recursive formulas as in \eqref{recuree1}, \eqref{recuree2} and \eqref{recuree3} for the solutions such that $\lambda(f)<\infty$, it is difficult to verify the existence of $b_2$ and $b_3$ satisfying these recursive formulas. We conjecture that equation \eqref{Bank-laine se3simple} can admit a nontrivial solution $f$ such that $\lambda(f)<\infty$ only when $l=2$ or $l=4$. We will study this conjecture further.

\section*{Acknowledgements}

The author is supported by a Project funded by China Postdoctoral Science Foundation~{(2020M680334)} and the Fundamental Research Funds for the Central Universities~{(FRF-TP-19-055A1)}. The author would like to thank professor Yik-man Chiang of the Hong Kong University of Science and Technology for sharing with the author the two references~\cite{Chiang2019} and~\cite{Kovacic1986}. This greatly simplifies the original proof of Theorem~\ref{maintheorem5}. The author also would like to thank the referee for his/her very valuable suggestions and comments.





\end{document}